\documentclass[11pt]{amsart}
\usepackage{amssymb,amsmath,a4wide,graphicx,stmaryrd,footmisc}
\usepackage[pagebackref=true]{hyperref}
\hypersetup{backref}
\usepackage[all]{xy}

\usepackage{color}

\title{Homotopy Theory of Labelled Symmetric Precubical Sets}

\author[P. Gaucher]{Philippe Gaucher}

\address{Laboratoire PPS  (CNRS UMR 7126)\\
Case 7014\\
Univ Paris Diderot\\
 Sorbonne Paris Cit\'e\\
 F-75205 PARIS\\ France}

\urladdr{http://www.pps.univ-paris-diderot.fr/{\~{}}gaucher/} 

\subjclass{18C35,18G55,55U35,68Q85}

\keywords{precubical set, higher dimensional transition system, locally presentable category, topological category, combinatorial model category, Bousfield localization}

\swapnumbers

\newcommand{\C}{\mathcal{C}}

\newcommand{\K}{\mathcal{K}}

\newcommand{\I}{\mathcal{I}}
\newcommand{\W}{\mathcal{W}}

\newcommand{\R}{\mathbb{R}}
\newcommand{\T}{\mathbb{T}}
\newcommand{\de}{\partial}
\newcommand{\p}\times

\newtheorem{thm}{Theorem}[section]
\newtheorem{prop}[thm]{Proposition}

\newtheorem{cor}[thm]{Corollary}

\newtheorem{defn}[thm]{Definition}
\newtheorem{nota}[thm]{Notation}
\newtheorem{rem}[thm]{Remark}
\newcommand{\bd}{\begin{defn}}
\newcommand{\ed}{\end{defn}}
\newcommand{\bp}{\begin{prop}}
\newcommand{\ep}{\end{prop}}
\newcommand{\bth}{\begin{thm}}
\renewcommand{\eth}{\end{thm}}
\newcommand{\bpf}{\begin{proof}}
\newcommand{\epf}{\end{proof}}
\newcommand{\bmul}{\begin{multline*}}
\newcommand{\emul}{\end{multline*}}
\newcommand{\bc}{\begin{cor}}
\newcommand{\ec}{\end{cor}}

\newcommand{\fL}[1]{\ar@{->}[ll]_-{#1}}
\newcommand{\fR}[1]{\ar@{->}[rr]^-{#1}}
\newcommand{\fRr}[1]{\ar@{->}[rrr]^-{#1}}
\newcommand{\fD}[1]{\ar@{->}[dd]_-{#1}}
\newcommand{\fU}[1]{\ar@{->}[uu]^-{#1}}
\newcommand{\f}[2]{\ar@{->}[#1]|{#2}}
\newcommand{\ff}[2]{\ar@2{->}[#1]|{#2}}
\newcommand{\frr}[1]{\ar@{->}[rrrr]^-{#1}}

\newcommand{\fl}[1]{\ar@{->}[l]_-{#1}}
\newcommand{\fr}[1]{\ar@{->}[r]^-{#1}}
\newcommand{\fd}[1]{\ar@{->}[d]_-{#1}}
\newcommand{\fu}[1]{\ar@{->}[u]^-{#1}}

\newcommand{\iso}{\cong}

\renewcommand{\leq}{\leqslant}
\renewcommand{\geq}{\geqslant}
\newcommand{\dd}[1]{\uparrow\!\!{#1}\!\!\uparrow}

\def\cartesien{%
  \ar@{-}[]+R+<6pt,-2pt>;[]+RD+<6pt,-6pt>%
  \ar@{-}[]+D+<2pt,-6pt>;[]+RD+<6pt,-6pt>%
}
\def\cocartesien{%
  \ar@{-}[]+L+<-6pt,+2pt>;[]+LU+<-6pt,+6pt>%
  \ar@{-}[]+U+<-2pt,+6pt>;[]+LU+<-6pt,+6pt>%
}
\def\hocartesien{%
  \ar@{-}[]+R+<6pt,-2pt>;[]+RD+<6pt,-6pt>_{h}%
  \ar@{-}[]+D+<2pt,-6pt>;[]+RD+<6pt,-6pt>%
}
\def\hococartesien{%
  \ar@{-}[]+L+<-6pt,+2pt>;[]+LU+<-6pt,+6pt>_{h}%
  \ar@{-}[]+U+<-2pt,+6pt>;[]+LU+<-6pt,+6pt>%
}

\newcommand{\brm}[1]{\rm{\mathbf{#1}}}

\newcommand{\set}{{\brm{Set}}}

\newcommand{\hda}{{\brm{HDA}}}

\DeclareMathOperator{\id}{Id}
\DeclareMathOperator{\sh}{Sh}
\DeclareMathOperator{\Mor}{Mor}

\newcommand{\liminj}{\varinjlim}

\DeclareMathOperator{\whdts}{\brm{WHDTS}}

\DeclareMathOperator{\cts}{\brm{CTS}}

\DeclareMathOperator{\cub}{\underline{Cub}}
\DeclareMathOperator{\cube}{CUBE}

\DeclareMathOperator{\bl}{\brm{\underline{L}}}
\newcommand{\ls}[1]{(#1/\!/\Sigma)}

\makeatletter
\def\varholim@#1#2{%
  \vtop{\m@th\ialign{##\cr
    \hfil$#1\operator@font holim$\hfil\cr
    \noalign{\nointerlineskip\kern1.5\ex@}#2\cr
    \noalign{\nointerlineskip\kern-\ex@}\cr}}%
}
\def\holimproj{%
  \mathop{\mathpalette\varholim@{\leftarrowfill@\textstyle}}\nmlimits@
}
\def\holiminj{%
  \mathop{\mathpalette\varholim@{\rightarrowfill@\textstyle}}\nmlimits@
}
\makeatother

\DeclareMathOperator{\cell}{{\brm{cell}}}
\DeclareMathOperator{\cof}{{\brm{cof}}}
\DeclareMathOperator{\proj}{{\brm{proj}}}
\DeclareMathOperator{\inj}{{\brm{inj}}}
\newcommand{\ddownarrow}{{\downarrow}}

\DeclareMathOperator{\cyl}{{Cyl}}

\DeclareMathOperator{\CSA}{CSA}

\begin{document}

\begin{abstract} 
  This paper is the third paper of a series devoted to higher
  dimensional transition systems. The preceding paper proved the
  existence of a left determined model structure on the category of
  cubical transition systems. In this sequel, it is proved that there
  exists a model category of labelled symmetric precubical sets which
  is Quillen equivalent to the Bousfield localization of this left
  determined model category by the cubification functor.  The
  realization functor from labelled symmetric precubical sets to
  cubical transition systems which was introduced in the first paper
  of this series is used to establish this Quillen equivalence.
  However, it is not a left Quillen functor. It is only a left
  adjoint. It is proved that the two model categories are related to
  each other by a zig-zag of Quillen equivalences of length two. The
  middle model category is still the model category of cubical
  transition systems, but with an additional family of generating
  cofibrations. The weak equivalences are closely related to
  bisimulation. Similar results are obtained by restricting the
  constructions to the labelled symmetric precubical sets satisfying
  the HDA paradigm.
\end{abstract}

\maketitle

\tableofcontents

\section{Introduction}

\subsection{Presentation of two combinatorial approaches of concurrency.}
\begin{figure}
\[
\xymatrix{
& () \ar@{->}[rd]^{(b)}&\\
()\ar@{->}[ru]^{(a)}\ar@{->}[rd]_{(b)} & (a,b) & ()\\
&()\ar@{->}[ru]_{(a)}&}
\]
\caption{Concurrent execution of two actions $a$ and $b$}
\label{concab}
\end{figure}
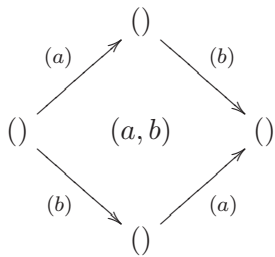

This paper is the third paper of a series devoted to higher
dimensional transition systems (HDTS). The first appearance of this
notion dates back to \cite{MR1461821} in which the concurrent
execution of $n$ action is modelled by a multiset of actions. The
first paper of this series \cite{hdts} proved that this approach is
actually the same as the geometric approach of concurrency dating back
to \cite{EWDCooperating} \cite{Pratt} \cite{Gunawardena1}.  It is
really not possible to give an exhaustive list of references for this
subject because this field of research is growing very fast but at
least two surveys are available \cite{survol} \cite{rvg} presenting
various topological and combinatorial models. In this theory, the
concurrent execution of $n$ actions is modelled by a labelled
$n$-cube. Each coordinate represents the state of one of the $n$
processes running concurrently, from $0$ (not started) to $1$
(complete). Figure~\ref{concab} represents the concurrent execution of
two actions $a$ and $b$.  

The formalism of labelled symmetric precubical sets is another example
of combinatorial object encoding these ideas. The idea of modelling
labelling \emph{cubical sets} by working in a comma category is
probably introduced in \cite{labelled}.  Labelled \emph{precubical}
sets, meaning without the standard degeneracy maps coming from
algebraic topology, and without symmetry operators, are actually
sufficient to model the geometry of all process algebras for any
synchronization algebra by \cite[Section~4]{ccsprecub}. Let us
emphasize this fact.  Not only are the standard degeneracy maps coming
from algebraic topology useless for modelling the space of paths of a
process algebra, but also new degeneracy maps, the \emph{transverse
  degeneracy maps} of \cite{symcub} seem to be required to better
understand the semantics of process algebras. These non-standard
degeneracy maps will not be used in this paper however. Every process
algebra can then be viewed as a labelled \emph{symmetric} precubical
set (Definition~\ref{def_lsps}) just by considering the \emph{free
  labelled symmetric precubical set} generated by the associated
labelled precubical set \cite{symcub}.  Thanks to the symmetry
operators, the parallel composition of two processes $P$ and $Q$ is
isomorphic to the parallel composition of two processes $Q$ and
$P$. Such an isomorphism just does not exist in general in the
category of labelled precubical sets, except in degenerate situations
like $P=Q$ of course.

A semantics of process algebras in terms of HDTS is expounded in
\cite{hdts}.  Unlike labelled symmetric precubical sets, HDTS do not
necessarily have face operators: they are not part of the definition
indeed (see Definition~\ref{weak-HDTS}). An immediate consequence is
that the colimit of the cubes contained in a given HDTS (called its
cubification, see Definition~\ref{def-cub}) does not necessarily give
back the HDTS. A nonempty HDTS may even have an empty
cubification. Even the cubical transition systems which are, by
definition, equal to the union of their subcubes are not necessarily
equal to their cubification (see below in this introduction).  There
is another striking difference between HDTS and labelled symmetric
precubical sets: all Cattani-Sassone higher dimensional transition
systems satisfy the so-called HDA paradigm (see Section~\ref{paradigm}
of this paper, and \cite[Section~7]{hdts}). This implies that the
formalization of the parallel composition, for any synchronization
algebra, of two processes is much simpler with HDTS than with
precubical sets, symmetric or not. Indeed, there is no need in the
setting of HDTS to introduce tricky combinatorial constructions like
the \emph{directed coskeleton construction} of
\cite{ccsprecub}~\footnote{Let us just recall here that the choice of
  ``directed'' in ``directed coskeleton construction'' was a very bad
  idea.}, or the \emph{transverse degeneracy maps} of
\cite{symcub}. We just have to list all higher dimensional transitions
of a parallel composition by reading the definition from a computer
science book and to put them in the set of transitions of the HDTS.

\subsection{The salient mathematical facts of the preceding papers of this series.} 
The first paper \cite{hdts} is devoted to introducing a more
convenient formalism to deal with HDTS. More precisely, the category
of weak HDTS is introduced (Definition~\ref{weak-HDTS}). It enjoys a
lot of very nice properties: topological~\footnote{A topological
  category is a category equipped with a forgetful topological functor
  towards a power of the category of sets.}, locally finitely
presentable. The category of Cattani-Sassone higher dimensional
transition systems is interpreted as a full reflective subcategory of
the category of weak HDTS \cite[Corollary~5.7]{hdts}. And it is proved
in \cite[Theorem~11.6]{hdts} that the categorical localization of the
category of Cattani-Sassone higher dimensional transition systems by
the cubification functor is equivalent to a full reflective
subcategory of that of labelled symmetric precubical sets. The main
tool is a realization functor from labelled symmetric precubical sets
to HDTS, whose construction is presented in an improved form in
Section~\ref{realization} of this paper.  Symmetry operators are
required for this result since the group of automorphisms of the
labelled $n$-cube in the category of Cattani-Sassone HDTS is the
$n$-th symmetry group, not the singleton as in the category of
labelled precubical sets of \cite{ccsprecub}. In other terms, the
category of HDTS has built-in symmetry operators which, of course,
come from the action of the $n$-th symmetric group on the set of
$n$-dimensional transitions. The realization functor from labelled
symmetric precubical sets to HDTS will be reused in this paper to get
Quillen equivalences.

The second paper of the series \cite{cubicalhdts} is devoted to the
study of the homotopy theory of HDTS. A left determined model
structure is built on the topological and finitely presentable
category of weak HDTS, and then restricted to the full subcategory of
cubical transition systems \cite[Corollary~6.8]{cubicalhdts}, i.e.
the weak HDTS which are equal to the union of their subcubes
(Definition~\ref{subcube}).  This full coreflective subcategory of
that of weak HDTS contains all examples coming from computer science
even if the topological structure of the larger category of weak HDTS
keeps playing an important role in the development of this theory. The
class of weak equivalences of this left determined model structure is
completely characterized. It appears that it is really too small to be
interesting. It also turns out that all weak equivalences are
bisimulations and it is tempting to Bousfield localize by all
bisimulations. By \cite[Theorem~9.5]{cubicalhdts}, such a Bousfield
localization exists but its study is out of reach at present. An
intermediate Bousfield localization, by the cubification functor
again, is proved to exist as well \cite[Section~8]{cubicalhdts}.

One word must be said about the notion of cubical transition
system. Not all HDTS are equal to the colimit of their subcubes. For
example the boundary $\de C_2[x_1,x_2]$ of the full $2$-cube
$C_2[x_1,x_2]$, which is obtained by removing all its $2$-transitions,
that is to say $((0,0),x_1,x_2,(1,1))$ and
$((0,0),x_2,x_1,(1,1))$~\footnote{The $n$-cube $C_n[x_1,\dots,x_n]$
  has actually and by definition $n$ distinct actions $(x_i,i)$ for
  $i=1,\dots ,n$, but it is assumed here that $x_1\neq x_2$, so there
  is no need to overload the notations by writing $(x_1,1)$,
  $(x_2,2)$; the $n$-cube viewed as a HDTS has exactly $n!$
  $n$-dimensional transitions.}. Indeed the HDTS $\de C_2[x_1,x_2]$
has only two actions $x_1$ and $x_2$, whereas its cubification has
four distinct actions $x_1,x_1',x_2,x'_2$ with the labelling map
$\mu(x_k)=\mu(x_k')=x_k$ for $k=1,2$. So $\de C_2[x_1,x_2]$ is not
isomorphic to its cubification, and therefore it cannot be a colimit
of cubes.  But it is cubical anyway. This is because of this subtle
point that we are forced to use the category of cubical transition
systems.  It is not known whether a model category structure like the
one of \cite{cubicalhdts} exists on the full subcategory of weak HDTS
of colimits of cubes. The main problem consists of finding another set
of generating cofibrations (instead of the set $\I$ defined in
Notation~\ref{cofgen}) without using the boundary of the labelled
$2$-cubes.

\subsection{Presentation of this paper.} 
This third paper of the series goes back to the link between labelled
symmetric precubical sets and cubical transition systems. One of the
main results of this paper is that a model category structure is
constructed on the category of labelled symmetric precubical sets
(Theorem~\ref{constr_model}) thanks to Marc Olschok's PhD \cite{MOPHD}
\cite{MO}.  And it is proved in Theorem~\ref{same1} that there exists
a Bousfield localization of the latter which is Quillen equivalent to
the model category structure of cubical transition systems introduced
in \cite{cubicalhdts} (not the left determined one, but its Bousfield
localization by the cubification functor, which is studied in
\cite[Section~8]{cubicalhdts}). It is remarkable that like for the
categorical equivalence of \cite{hdts}, the cubification functor is
used once again, this time to obtain a Quillen equivalence.
Theorem~\ref{same2} is a similar theorem after restriction to the
labelled symmetric precubical sets satisfying the HDA paradigm.
Unfortunately, almost nothing is known about these Bousfield
localizations.

Surprisingly, the realization functor from labelled symmetric
precubical sets to cubical transition systems is not a left Quillen
functor. It is only a left adjoint
(Proposition~\ref{not-left-Quillen}). An intermediate model category
must be used in the proofs to get the Quillen equivalences
(Theorem~\ref{homotopy-of-cts-plus}). The cause of this problem is the
family of cofibrations consisting of the inclusions $\de\square_S[x,y]
\subset \square_S[x,y]$ of the boundary of a labelled $2$-cube to the
full labelled $2$-cube (Proposition~\ref{not-left-Quillen}) for $x$
and $y$ running over the set of labels $\Sigma$. The image by the
realization functor is not a cofibration of HDTS. Indeed, the image of
$\de\square_S[x,y]$ is precisely the cubification of the boundary of
the $2$-cube $\de C_2[x,y]$ because every labelled symmetric
precubical set is equal to the colimit of its cubes and because the
realization functor is colimit-preserving. It has four actions (see
above !)  $x,x',y,y'$ with the labelling map $\mu(x)=\mu(x')=x$ and
$\mu(y)=\mu(y')=y$ whereas the realization of $\square_S[x,y]$ is the
$2$-cube $C_2[x,y]$ which has two actions $x$ and $y$: therefore the
map from the realization of $\de\square_S[x,y]$ to the one of
$\square_S[x,y]$ cannot be one-to-one on actions, so it cannot be a
cofibration of HDTS by definition. The intermediate model category is
precisely obtained by adding this family of maps (so the realization
of the inclusions $\de\square_S[x,y] \subset \square_S[x,y]$) to the
set of generating cofibrations of the model category of cubical
transition systems ! In other terms, we force the realization functor
from labelled symmetric precubical sets to cubical transition systems
to become a left Quillen functor. And the second surprise is that that
just works fine.

Again the same family of inclusions $\de\square_S[x,y] \subset
\square_S[x,y]$ for $x$ and $y$ running over the set of labels
$\Sigma$ prevents the interval object of labelled symmetric precubical
sets from being very good. It is only good (Proposition~\ref{wfs_OK}
and the remark after the proof). The realization as HDTS of the same
family of cofibrations also prevents the interval object of cubical
transition systems from being very good as well with respect to the
augmented set of generating cofibrations
(Theorem~\ref{homotopy-of-cts-plus}). As a consequence, the Olschok
construction cannot tell us anything about the left determinedness of
the model category of labelled symmetric precubical sets and of the
augmented model category of cubical transition systems.

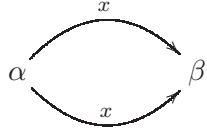
\begin{figure}
\[
\xymatrix{
\alpha \ar@/^20pt/[rr]^-{x} \ar@/_20pt/[rr]^-{x}&& \beta }
\]
\caption{$\cyl(\square_S[x])$: the cylinder of $\square_S[x]$ is
  homotopy equivalent to the $1$-cube $\square_S[x]$.}
\label{Da0}
\end{figure}

This new model category structure on labelled symmetric precubical
sets is very different from the ones coming from algebraic
topology. Indeed, the class of cofibrations is strictly larger than
the class of monomorphisms. Like the model category of flows
\cite{model3}, it contains the generating cofibration $R:\{0,1\} \to
\{0\}$. This makes it impossible to use tools like Cisinski's homotopy
theory of toposes \cite{MR1924082} or Hirschhorn's theory of Bousfield
localization \cite{ref_model2}. The main technical tool of this paper
is Marc Olschok's PhD thesis \cite{MOPHD} \cite{MO} instead. Moreover,
not only is the $1$-cube not weakly equivalent to a point; it is in
fact weakly equivalent to two copies of itself where the two initial
(final resp.) states are identified as in Figure~\ref{Da0}. This new
model category is adapted, like the ones constructed on cubical
transition systems in \cite{cubicalhdts}, to the study of bisimulation
\cite{MR1365754} \cite{0856.68067}. In the case of Figure~\ref{Da0},
the labelled symmetric precubical set has the same behavior as the
$1$-cube $\square_S[x]$ labelled by $x$. Indeed, the unique map from
$\cyl(\square_S[x])$ to $C_1[x]$ is a bisimulation.

\subsection*{Outline of the paper}

Section~\ref{Olschok} is a reminder about the Olschok construction of
combinatorial model categories, at least the first part of his PhD
devoted to the generalization of Cisinski's work to the setting of
locally presentable categories. Only what is used in this paper is
recalled. So the statement of Theorem~\ref{build-model-cat} is
certainly less general than what is written in \cite{MO} and
\cite{MOPHD}.  Section~\ref{reminder-cts} is a reminder about weak
HDTS and cubical transition systems. Several important basic examples
of such objects are given.  Section~\ref{homotopy-cts} recalls the
homotopy theory of cubical transition systems. The exposition is
improved, so it is more than just a reminder. In particular, an
explicit set of generating cofibrations is given.
Section~\ref{rappel_sym_pre_set} recalls the definition of labelled
symmetric precubical set. Once again, several important basic examples
are given.  Section~\ref{homotopy-cube} constructs the new model
category structure on labelled symmetric precubical sets
(Theorem~\ref{constr_model}). Roughly speaking, we \emph{really} just
have to mimic the construction of the model category structure on
cubical transition systems. Section~\ref{realization} recalls the
construction of the realization functor from labelled symmetric
precubical sets to cubical transition systems. The exposition is much
better than in \cite{hdts} where it is introduced, so it is also more
than just a reminder. It is also proved in the same section that the
realization functor is not a left Quillen functor, and it is explained
how to overcome this problem by adding one family of generating
cofibrations to the category of cubical transition systems.  And
Section~\ref{homotopy-realization} proves one of the main result of
this paper: there exists a model category of labelled symmetric
precubical sets which is Quillen equivalent to the homotopy theory of
cubical transition systems (Theorem~\ref{same1}). The last section
restricts the homotopy constructions to the full reflective
subcategory of labelled symmetric precubical sets satisfying the HDA
paradigm (Theorem~\ref{constr_model_HDA}) and proves a similar result
(Theorem~\ref{same2}).

\subsection*{Preliminaries}

The necessary bibliographical references and reminders are given
throughout the text. The category of sets is denoted by $\set$.  All
categories are locally small. The set of maps in a category $\K$ from
$X$ to $Y$ is denoted by $\K(X,Y)$. The locally small category whose
objects are the maps of $\K$ and whose morphisms are the commutative
squares is denoted by $\Mor(\K)$. The initial (final resp.) object, if
it exists, is always denoted by $\varnothing$ ($\mathbf{1}$). The
identity of an object $X$ is denoted by $\id_X$.  A subcategory is
always \emph{isomorphism-closed}.  We refer to \cite{MR95j:18001} for
locally presentable categories, to \cite{MR2506258} for combinatorial
model categories, and to \cite{topologicalcat} for topological
categories (i.e. categories equipped with a topological functor
towards a power of the category of sets).  We refer to
\cite{MR99h:55031} and to \cite{ref_model2} for model categories. For
general facts about weak factorization systems, see also
\cite{ideeloc}. We recommend the reading of Marc Olschok's PhD
\cite{MOPHD}. The first part, published in \cite{MO}, is used in this
paper.

\subsection*{Acknowledgements}
The author would like to thank the referee for suggestions on
improving the exposition.

\section{The Olschok construction of model categories}
\label{Olschok}

We want to review the Olschok construction of combinatorial model
categories \cite{MO}, as it is already used in \cite{cubicalhdts}, and
as it is used in this paper, i.e. by starting from a \emph{good
  interval object}, i.e. a good cylinder functor $\cyl(-)$ of the form
$\cyl(-)=V\p -$ where $V$ is a bipointed object of the ambient
category.

Let $f$ and $g$ be two maps of a locally presentable category
$\K$. Write $f\square g$ when $f$ satisfies the \emph{left lifting
  property} with respect to $g$ (or equivalently $g$ satisfies the
\emph{right lifting property} with respect to $f$). Let us introduce
the notations $\inj_\K(\C) = \{g \in \K, \forall f \in \C, f\square
g\}$, $\proj_\K(\C) = \{f \in \K, \forall g \in \C, f\square g\}$ and
$\cof_\K(\C) = \proj_\K(\inj_\K(\C))$ where $\C$ is a class of maps of
a locally presentable category $\K$. The class of morphisms of $\K$
that are transfinite compositions of pushouts of elements of $\C$ is
denoted by $\cell_\K(\C)$.

\begin{nota} For every map $f:X \to Y$ and every natural
  transformation $\alpha : F \to F'$ between two endofunctors
  of $\K$, the map $f\star \alpha$ is the canonical map
\[f\star \alpha : FY \sqcup_{FX} F'X \longrightarrow F'Y\] 
induced by the commutative diagram of solid arrows 
\[\xymatrix{
FX \fR{\alpha_X}\fD{Ff} && F'X \fD{F'f} \\
&& \\
FY \fR{\alpha_Y} && F'Y}
\]
and the universal property of the pushout. For a set of morphisms
$\mathcal{A}$, let $\mathcal{A} \star \alpha = \{f\star \alpha, f\in
\mathcal{A}\}$.
\end{nota}

\bd Let $I$ be a set of maps of a locally presentable category $\K$.
A {\rm good cylinder} with respect to $I$ is a functor $\cyl:\K\to\K$
together with two natural transformations $\gamma^k:\id\Rightarrow
\cyl$ for $k=0,1$ and a natural transformation $\sigma : \cyl
\Rightarrow \id$ such that the codiagonal $\id\sqcup \id \Rightarrow \id$ 
factors as a composite 
\[
\xymatrix
{
\id\sqcup \id\ar@{=>}[rr]^-{\gamma^0\sqcup \gamma^1}&& \cyl \ar@{=>}[rr]^-{\sigma}&& \id
}
\]
and such that the left-hand natural transformation $\gamma^0\sqcup
\gamma^1$ induces for all $X\in \K$ a map
\[\xymatrix{X\sqcup X \fR{\gamma^0_X \sqcup \gamma^1_X} &&  \cyl(X) \in \cof_\K(I).}\] When
moreover the right-hand map $\sigma_X$ belongs to $\inj_\K(I)$ for all
$X\in \K$, the functor $\cyl$ is called a {\rm very good
  cylinder}. \ed

\bd Let $I$ be a set of maps of a locally presentable category $\K$. A
good cylinder $\cyl:\K \to \K$ with respect to $I$ is {\rm cartesian}
if it is exponential  and if there are the inclusions $\cof_\K(I) \star
\gamma \subset \cof_\K(I)$ and $\cof_\K(I) \star \gamma^k \subset
\cof_\K(I)$ for $k=0,1$ where $\gamma = \gamma^0 \sqcup \gamma^1$.
\ed

In this paper, all cylinders will be of the form the binary product by
a bipointed object $\gamma^0,\gamma^1:\mathbf{1} \rightrightarrows V$
called the \emph{interval object}.  The natural transformations
$\gamma^k:\id \Rightarrow \cyl$ are equal to the natural
transformations $(\mathbf{1}\p -) \Rightarrow (V\p -)$ induced by the
two maps $\gamma^k : \mathbf{1} \to V$ for $k=0,1$. The natural
transformation $ \sigma:\cyl \Rightarrow \id$ is equal to the natural
transformation $(V\p -)\Rightarrow (\mathbf{1}\p -)$ induced by the
map $\sigma:V\to \mathbf{1}$ which is the unique map from the interval
object to the terminal object.  An interval object will be good (very
good, cartesian resp.) if and only if the corresponding cylinder
functor is good (very good, cartesian resp.). 

\begin{nota} Let $I$ and $S$ be two sets of maps of a locally
  presentable category $\K$. Let $V$ be a good interval object
  with respect to $I$. Denote sets of maps $\Lambda_\K^n(V,S,I)$
  recursively by $\Lambda^0_\K(V,S,I) = S \cup (I \star \gamma^0)
  \cup (I \star \gamma^1)$ and $\Lambda^{n+1}_\K(V,S,I) =
  \Lambda^n_\K(V,S,I) \star \gamma$ for $n \geq 0$. Then let
  \[\Lambda_\K(V,S,I) = \bigcup_{n\geq 0} \Lambda^n_\K(V,S,I).\]
\end{nota}

Let $\simeq$ be the homotopy relation associated with the cylinder
$V\p -$, i.e. for all maps $f,g:X\to Y$, $f\simeq g$ is equivalent to
the existence of a \emph{homotopy} $H:V\p X \to Y$ with $H\circ
\gamma^0 = f$ and $H\circ \gamma^1 = g$.

\begin{nota} We denote by $\W_\K(V,S,I)$ the class of maps $f : X \to
  Y$ of $\K$ such that for every $\Lambda_\K(V,S,I)$-injective object
  $T$, the induced set map \[\K(Y,T)/\!\!\simeq
  \stackrel{\iso}\longrightarrow \K(X,T)/\!\!\simeq\] is a bijection.
\end{nota}

We are now ready to recall the Olschok construction for this
particular setting:

\bth \label{build-model-cat} (Olschok) Let $\K$ be a locally
presentable category.  Let $I$ be a set of maps of $\mathcal{K}$. Let
$S \subset \cof_\K(I)$ be an arbitrary set of maps of $\K$. Let $V$ be
a good cartesian interval object with respect to $I$.  Suppose also
that for any object $X$ of $\K$, the canonical map $\varnothing \to X$
belongs to $\cof_\K(I)$. Then there exists a unique combinatorial
model category structure with class of cofibrations $\cof_\K(I)$ such
that the fibrant objects are the $\Lambda_\K(V,S,I)$-injective
objects. The class of weak equivalences is $\W_\K(V,S,I)$. All objects
are cofibrant.  \eth

\bpf Since all objects are cofibrant, the class of weak equivalences
is necessarily $\W_\K(V,S,I)$ by
\cite[Theorem~7.8.6]{ref_model2}. Hence the uniqueness. The existence
is a consequence of \cite[Theorem~3.16]{MO}.  \epf

\begin{nota} For $S=\varnothing$, the model category is just denoted by 
$\K$. \end{nota}

If the interval is very good in Theorem~\ref{build-model-cat}, then
$\W_\K(V,S,I)$ is the localizer generated by $S$ (with respect to the
class of cofibrations $\cof_\K(I)$) by \cite[Theorem~4.5]{MO} and $\K$
is then left determined in the sense of \cite{rotho}.  And the model
category we obtain for $S\neq\varnothing$ is the Bousfield
localization $\bl_{\mathcal{S}}(\K)$ of the left determined one by the
set of maps $S$. If the interval is only good, then the Olschok
construction can only tell us that the model category we obtained is
the Bousfield localization $\bl_{\Lambda_\K(V,S,I)}(\K)$ with respect
to $\Lambda_\K(V,S,I)$ because, by \cite[Lemma~4.4]{MO}, the class of
maps $\W_\K(V,S,I)$ is the localizer generated by $\Lambda_\K(V,S,I)$.

\section{Cubical transition systems}
\label{reminder-cts}

\begin{nota} A nonempty set of {\rm labels} $\Sigma$ is fixed.  \end{nota}

\bd \label{weak-HDTS} A {\rm weak higher dimensional transition system
  (weak HDTS)} consists of a triple \[(S,\mu:L\to
\Sigma,T=\bigcup_{n\geq 1}T_n)\] where $S$ is a set of {\rm states},
where $L$ is a set of {\rm actions}, where $\mu:L\to \Sigma$ is a set
map called the {\rm labelling map}, and finally where $T_n\subset S\p
L^n\p S$ for $n \geq 1$ is a set of {\rm $n$-transitions} or {\rm
  $n$-dimensional transitions} such that one has:
\begin{itemize}
\item (Multiset axiom) For every permutation $\sigma$ of
  $\{1,\dots,n\}$ with $n\geq 2$, if $(\alpha,u_1,\dots,u_n,\beta)$ is
  a transition, then $(\alpha,u_{\sigma(1)}, \dots, u_{\sigma(n)},
  \beta)$ is a transition as well.
\item (Composition axiom) For every $(n+2)$-tuple
  $(\alpha,u_1,\dots,u_n,\beta)$ with $n\geq 3$, for every $p,q\geq 1$
  with $p+q<n$, if the five tuples 
\begin{multline*}
  (\alpha,u_1, \dots, u_n, \beta), (\alpha,u_1, \dots, u_p, \nu_1), 
  (\nu_1, u_{p+1}, \dots, u_n, \beta),\\ (\alpha, u_1, \dots, u_{p+q},
  \nu_2),  (\nu_2, u_{p+q+1}, \dots, u_n, \beta)
\end{multline*}
are transitions, then the $(q+2)$-tuple $(\nu_1, u_{p+1}, \dots,
u_{p+q}, \nu_2)$ is a transition as well.
\end{itemize}
A map of weak higher dimensional transition systems
\[f:(S,\mu : L \to \Sigma,(T_n)_{n\geq 1}) \to
(S',\mu' : L' \to \Sigma ,(T'_n)_{n\geq 1})\] consists of a
set map $f_0: S \to S'$, a commutative square
\[
\xymatrix{
  L \ar@{->}[r]^-{\mu} \ar@{->}[d]_-{\widetilde{f}}& \Sigma \ar@{=}[d]\\
  L' \ar@{->}[r]_-{\mu'} & \Sigma}
\] 
such that if $(\alpha,u_1,\dots,u_n,\beta)$ is a transition, then
$(f_0(\alpha),\widetilde{f}(u_1),\dots,\widetilde{f}(u_n),f_0(\beta))$
is a transition. The corresponding category is denoted by $\whdts$.
The $n$-transition $(\alpha,u_1,\dots,u_n,\beta)$ is also called a
{\rm transition from $\alpha$ to $\beta$}.  \ed

\begin{nota} The labelling map from the set of actions to the set of
  labels will be very often denoted by $\mu$. \end{nota}

The category $\whdts$ is locally finitely presentable by
\cite[Theorem~3.4]{hdts}. The functor \[\omega : \whdts
\longrightarrow \set^{\{s\}\cup \Sigma}\] taking the weak higher
dimensional transition system $(S,\mu : L \to \Sigma,(T_n)_{n\geq 1})$
to the $(\{s\}\cup \Sigma)$-tuple of sets $(S,(\mu^{-1}(x))_{x\in
  \Sigma}) \in \set^{\{s\}\cup \Sigma}$ is topological by
\cite[Theorem~3.4]{hdts} too.

There is a slight change in the terminology with respect to the one of
\cite{hdts} and \cite{cubicalhdts}. The \emph{Coherence axiom} is
called now the \emph{Composition axiom} because this axiom really
looks like a $5$-ary composition even if it is not known what
conclusion should be drawn from such an observation.

\begin{nota} For $n\geq 1$, let $0_n = (0,\dots,0)$ ($n$-times) and
  $1_n = (1,\dots,1)$ ($n$-times). By convention, let 
  $0_0=1_0=()$. \end{nota}

We give now some important examples of weak HDTS. In each of the
following examples, the Multiset axiom and the Composition axiom are
satisfied for trivial reasons.

\begin{enumerate}
\item Let $n \geq 0$. Let $x_1,\dots,x_n \in \Sigma$. The \emph{pure
    $n$-transition} $C_n[x_1,\dots,x_n]^{ext}$ is the weak HDTS with
  the set of states $\{0_n,1_n\}$, with the set of actions $\{(x_1,1),
  \dots, (x_n,n)\}$ and with the transitions all $(n+2)$-tuples
  \[(0_n,(x_{\sigma(1)},\sigma(1)), \dots,
  (x_{\sigma(n)},\sigma(n)),1_n)\] for $\sigma$ running over the set of
  permutations of the set $\{1,\dots ,n\}$.
\item Every set $X$ may be identified with the weak HDTS having the
  set of states $X$, with no actions and no transitions.
\item For every $x\in \Sigma$, let us denote by $\underline{x}$ the
  weak HDTS with no states, one action $x$, and no
  transitions. Warning: the weak HDTS $\{x\}$ contains one state $x$
  and no actions whereas the weak HDTS $\underline{x}$ contains no
  states and one action $x$.
\end{enumerate}

The following example plays a special role in the theory: 

\begin{nota} \label{doublet} For every $x\in \Sigma$, let us denote by
$\dd{x}$ the weak HDTS with four states $\{1,2,3,4\}$, one action $x$
and two transitions $(1,x,2)$ and $(3,x,4)$.
\end{nota}

Another important example is the one of the $n$-cube which is recalled
now.

\bp \label{cas_cube} \cite[Proposition~5.2]{hdts} Let $n\geq 0$ and
$x_1,\dots,x_n\in \Sigma$. Let $T_d\subset \{0,1\}^n \p
\{(x_1,1),\dots,(x_n,n)\}^d \p \{0,1\}^n$ (with $d\geq 1$) be the
subset of $(d+2)$-tuples
\[((\epsilon_1,\dots,\epsilon_n), (x_{i_1},i_1),\dots,(x_{i_d},i_d),
(\epsilon'_1,\dots,\epsilon'_n))\] such that
\begin{itemize}
\item $i_m = i_n$ implies $m = n$, i.e. there are no repetitions in the
  list \[(x_{i_1},i_1),\dots,(x_{i_d},i_d)\]
\item for all $i$, $\epsilon_i\leq \epsilon'_i$
\item $\epsilon_i\neq \epsilon'_i$ if and only if
  $i\in\{i_1,\dots,i_d\}$. 
\end{itemize}
Let $\mu : \{(x_1,1),\dots,(x_n,n)\} \to \Sigma$ be the set
map defined by $\mu(x_i,i) = x_i$. Then \[C_n[x_1,\dots,x_n] =
(\{0,1\}^n,\mu : \{(x_1,1),\dots,(x_n,n)\}\to
\Sigma,(T_d)_{d\geq 1})\] is a well-defined weak HDTS called the {\rm
  $n$-cube}. \ep

For $n = 0$, $C_0[]$, also denoted by $C_0$, is nothing else but the
weak HDTS $(\{()\},\mu:\varnothing \to
\Sigma,\varnothing)$. For every $x\in \Sigma$, one has $C_1[x] =
C_1[x]^{ext}$.

\bd \label{boundary-def} Let $n\geq 1$ and $x_1,\dots,x_n \in
\Sigma$. Let $\de C_n[x_1,\dots,x_n]$ be the weak HDTS defined by
removing from its set of transitions all $n$-transitions. It is called
the {\rm boundary} of $C_n[x_1,\dots,x_n]$. \ed

We restricted our attention in \cite{cubicalhdts} to the so-called
cubical transition systems, i.e. the weak HDTS which are equal to the
union of their subcubes. These weak HDTS include all useful examples.

\bd \label{subcube} Let $X$ be a weak HDTS. A {\rm cube} of $X$ is a
map $C_n[x_1,\dots,x_n] \longrightarrow X$. A {\rm subcube} of $X$ is
the image of a cube of $X$. A weak HDTS is a {\rm cubical transition
  system} if it is equal to the union of its subcubes. The full
subcategory of cubical transition systems is denoted by $\cts$.  \ed

Note that the weak HDTS $\de C_2[x_1,x_2]$ is not a colimit of cubes
but is cubical (see \cite[Corollary~3.12]{cubicalhdts} and the
discussion after it): it is obtained by identifying states in the
cubical transition system $\dd{x_1} \sqcup \dd{x_2}$. This is the
reason why we do not work in \cite{cubicalhdts} with the subcategory
of \emph{colimits} of cubes.

The category $\cts$ is a small-injectivity class, and a full
\emph{coreflective} locally finitely presentable subcategory of
$\whdts$ by \cite[Corollary~3.15]{cubicalhdts}. More precisely, a weak
HDTS is cubical if and only if it is injective with respect to the
maps $\underline{x} \subset C_1[x]$ for all $x\in \Sigma$ and to the
maps $C_n[x_1,\dots,x_n]^{ext} \subset C_n[x_1,\dots,x_n]$ for all $n
\geq 0$ and $x_1,\dots,x_n\in \Sigma$ by
\cite[Theorem~3.6]{cubicalhdts}.

\bd Let $X$ be a weak HDTS. An action $u$ of $X$ is {\rm used} if
there exists a $1$-transition $(\alpha,u,\beta)$. {\rm All actions are
  used} if $X$ is injective with respect to the maps $\underline{x}
\subset C_1[x]$ for all $x\in \Sigma$. \ed

\bd \label{isa} A weak HDTS $X$ satisfies the {\rm Intermediate state
  axiom} if for every $n\geq 2$, every $p$ with $1\leq p<n$ and every
transition $(\alpha,u_1,\dots,u_n,\beta)$ of $X$, there exists a (not
necessarily unique) state $\nu$ such that both
$(\alpha,u_1,\dots,u_p,\nu)$ and $(\nu,u_{p+1},\dots,u_n,\beta)$ are
transitions.  \ed

By \cite[Proposition~6.6]{cubicalhdts}, a weak HDTS satisfies the
Intermediate state axiom if and only if it is injective with respect
to the maps $C_n[x_1,\dots,x_n]^{ext} \subset C_n[x_1,\dots,x_n]$ for
all $n \geq 0$ and $x_1,\dots,x_n\in \Sigma$.  So a weak HDTS is
cubical if and only if all actions are used and it satisfies the
Intermediate state axiom.

\section{The homotopy theory of cubical transition systems}
\label{homotopy-cts}

Let us recall now the homotopy theory of $\cts$. This third paper
about higher dimensional transition systems contains some improvements
in the exposition of this theory. In particular, a set of generating
cofibrations can now be exhibited (in \cite{cubicalhdts}, the
existence of a set of generating cofibrations is proved using
transfinite techniques).

\bd A {\rm cofibration} of cubical transition systems is a map of weak
HDTS inducing an injection between the set of actions.  \ed

To make the reading of this paper easier, let us introduce a new
notation (which will be used later in Proposition~\ref{p0}).

\begin{nota} Let $\mathbb{S} : \set\ddownarrow \Sigma \to \cts$ be the
  functor given on objects as follows: if $\mu:L \to \Sigma$ is a set
  map then $\mathbb{S}(\mu)$ is the weak HDTS with set of states
  $\{0\}$, with labelling map $\mu$, and with set of
  transitions \[\{0\} \p \bigcup_{n\geq 1} L^n \p \{0\}.\]
\end{nota}

Note here that $\mathbb{S}(\mu)$ is a cubical transition system
because all actions are used and the Intermediate state axiom is
satisfied.

\bd Let us call $V:= \mathbb{S}(\Sigma\p\{0,1\} \to \Sigma)$ the {\rm
  interval object} of $\cts$ where $\Sigma\p\{0,1\} \to \Sigma$ is the
projection map. \ed

Note that $\mathbb{S}(\id_\Sigma)$ is the terminal object $\mathbf{1}$
of $\cts$.  For $k \in \{0,1\}$, denote by $\gamma^k:\mathbf{1} \to V$
the map of cubical transition systems induced by the composite set map
$\Sigma \iso \Sigma \p \{k\} \subset \Sigma \p \{0,1\}$. And denote by
$\sigma:V\to \mathbf{1}$ the canonical map, also induced by the
projection map $\Sigma \p \{0,1\} \to \Sigma$. Let $\gamma = \gamma^0
\sqcup \gamma^1 : \mathbf{1} \sqcup \mathbf{1} \to V$.

The interval $V$ is exponential by
\cite[Proposition~5.8]{cubicalhdts}. It is very good by
\cite[Proposition~5.7]{cubicalhdts} and cartesian by
\cite[Proposition~5.10]{cubicalhdts}. We are going to use the
following fact which is already implicitly present in \cite{hdts} and
\cite{cubicalhdts}.

\bp \label{transition-inj} Let $f: A\to B$ be a map of weak HDTS which
is bijective on states and actions.  Then it is one-to-one on
transitions. \ep

\bpf Let $(\alpha,u_1,\dots,u_m,\beta)$ and
$(\alpha',u'_1,\dots,u'_{m'},\beta')$ be two transitions of $A$ such
that
\[(f_0(\alpha),\widetilde{f}(u_1),\dots,\widetilde{f}(u_m),f_0(\beta))
=
(f_0(\alpha'),\widetilde{f}(u'_1),\dots,\widetilde{f}(u'_{m'}),f_0(\beta')).\] Then
$m = m'$, $f_0(\alpha) = f_0(\alpha')$, $\widetilde{f}(u_i) =
\widetilde{f}(u'_i)$ for $1\leq i \leq n$ and $f_0(\beta) =
f_0(\beta')$. So by hypothesis, $\alpha = \alpha'$, $\beta = \beta'$
and $u_i = u'_i$ for $1\leq i \leq n$. Hence
\[(\alpha,u_1,\dots,u_m,\beta) = (\alpha',u'_1,\dots,u'_{m'},\beta').\]
\epf

\begin{nota} \label{cofgen} Let $\I$~\footnote{The notations
    $C:\varnothing \to \{0\}$ and $R:\{0,1\} \to \{0\}$ are already
    used in \cite{model3} and in \cite{interpretation-glob} for the
    same generating cofibrations (in different categories of
    course). We will stick to this notation here, and for the model
    category of labelled symmetric precubical sets as well.} be the
  set of maps of cubical transition systems given by
\begin{multline*}
\I = \{C:\varnothing \to \{0\}, R:\{0,1\} \to \{0\}\} \\ \cup \{\de C_n[x_1,\dots,x_n]
\to C_n[x_1,\dots,x_n]\mid \hbox{$n\geq 1$ and $x_1,\dots,x_n \in
  \Sigma$}\} \\ \cup \{C_1[x] \to \dd{x}\mid x\in \Sigma\}.
\end{multline*}
\end{nota}

Let us recall again that the cubical HDTS $ \dd{x}$ is not a colimit
of cubes by \cite[After Definition~3.13 and before
  Remark~3.14]{cubicalhdts}. The colimit of all cubes (left-hand
cubical transition system of Figure~\ref{contrex-mono}) of $\dd{x}$
(right-hand cubical transition system of Figure~\ref{contrex-mono}) is
equal to the coproduct of two copies of $C_1[x]$ and it has two
actions $x_1$ and $x_2$ with $\mu(x_1)=\mu(x_2)=x$ whereas $\dd{x}$
has only one action $x$. Figure~\ref{contrex-mono} is an example of a
map of $\cts$ which is not an isomorphism, but which is a mono and an
epi.  Hence the additional family of cofibrations $C_1[x] \to \dd{x}$
for $x$ running over $\Sigma$ cannot be deduced from the other
generating cofibrations.

\begin{figure}
\[\left\{\begin{array}{c} C_1[x] \sqcup C_1[x] \\
        {\stackrel{x_1}\longrightarrow} \\
        {\stackrel{x_2}\longrightarrow}\end{array}\right.
  \stackrel{p_{x}}\longrightarrow \left\{\begin{array}{c} \liminj (C_1[x] \leftarrow \underline{x}
        \rightarrow C_1[x]) \\ {\stackrel{x}\longrightarrow}\\
        {\stackrel{x}\longrightarrow}\end{array}\right.\]
\caption{Monomorphism and epimorphism in $\cts$ with $\mu(x_1) =
  \mu(x_2) = x$}
\label{contrex-mono}
\end{figure}

\bth \label{generators} The class of cofibrations of cubical
transition systems is equal to \[\cell_{\cts}(\I) =
\cof_{\cts}(\I).\] \eth

\bpf Let $f:A \to B$ be a cofibration of cubical transition systems,
i.e. a map of cubical transition systems which is one-to-one on
actions. Let us factor $f$ as a composite $A \to Z \to B$ where the
left-hand map belongs to $\cell_{\cts}(\{C,R\})$ and where the
right-hand map belongs to $\inj_{\cts}(\{C,R\})$. Then the sets of
states of $Z$ and $B$ coincide, therefore we can suppose without lack
of generality that $f$ induces a bijection between the sets of
states. For every action $u$ of $B$ which does not belong to $A$, the
map $\underline{u} \to B$ factors (not in a unique way) as a composite
$\underline{\mu(u)} \to C_1[\mu(u)] \to B$ because $B$ is cubical.
The map $f:A \to B$ is bijective on states, therefore the composite
$C_0 \sqcup C_0 \iso \de C_1[\mu(u)] \subset C_1[\mu(u)] \to B$
factors as a composite $C_0 \sqcup C_0 \to A \to B$. Then for every
action $u$ of $B$ not in $A$, there exists a commutative square
\[
\xymatrix{C_0 \sqcup C_0 \iso \de C_1[\mu(u)] \fD{f_u} \fR{g_u} && A \fD{} \\
&& \\
  C_1[\mu(u)] \fR{} && B.
}
\]
Then consider the pushout diagram
\[
\xymatrix{\bigsqcup_{u}\de C_1[\mu(u)] \fD{f_u} \fR{g_u} && A \fD{} \\
&& \\
 \bigsqcup_{u} C_1[\mu(u)] \fR{} && \cocartesien Z.
}
\]
The canonical map $Z\to B$ induced by the pushout is bijective
both on states and on actions, and by
Proposition~\ref{transition-inj}, injective on transitions. Let us now
factor the map $Z \to B$ as a composite $Z \to D
\to B$ where the left-hand map belongs to 
\begin{multline*}\cell_{\cts}(\left\{\de C_n[x_1,\dots,x_n] \to
C_n[x_1,\dots,x_n]\mid n\geq 2 \hbox{ and }x_1,\dots,x_n \in \Sigma\right\} \\ 
\cup \left\{C_0 \sqcup C_0 \sqcup C_1[x] \to \dd{x}\mid x\in \Sigma\right\}) \end{multline*}
and the right-hand map belongs to 
\begin{multline*}\inj_{\cts}(\left\{\de C_n[x_1,\dots,x_n] \to
    C_n[x_1,\dots,x_n]\mid n\geq 2 \hbox{ and }x_1,\dots,x_n \in
    \Sigma\right\} \\\cup \left\{C_0 \sqcup C_0 \sqcup C_1[x]
    \to \dd{x}\mid x\in \Sigma\right\}).\end{multline*} where
the map $C_0 \sqcup C_0 \sqcup C_1[x] \to \dd{x}$ is defined
so that it is bijective on states. It is important to notice that the
maps $\de C_n[x_1,\dots,x_n] \to C_n[x_1,\dots,x_n]$ for every
$n\geq 2$ and every $x_1,\dots,x_n \in \Sigma$ and the maps $C_0
\sqcup C_0 \sqcup C_1[x] \to \dd{x}$ for every $x\in \Sigma$
are bijective on states and actions. Therefore the two maps of cubical
transition systems $Z \to D$ and $D \to B$ are
bijective on states and actions, and by
Proposition~\ref{transition-inj} injective on transitions.

Let $(\alpha,u,\beta)$ be a $1$-transition of $B$. Then $u$ is an
action of $B$ and therefore of $D$ as well and $\alpha$ and $\beta$
are two states of $D$. Since $D$ is cubical, there exists a
$1$-transition $(\alpha',u,\beta')$ of $D$ giving rise to a map
$C_1[\mu(u)] \to D$. Then consider the commutative diagram
\[
\xymatrix{\{\alpha\} \sqcup \{\beta\} \sqcup C_1[\mu(u)] \fD{} \fR{} && D \fD{} \\
&& \\
 \dd{\mu(u)} \fR{} \ar@{-->}[rruu]^{k} && B.
}
\]
The existence of the lift $k$ implies that the transition
$(\alpha,u,\beta)$ belongs to $D$, hence the map $D \to B$ is
onto on $1$-transitions. Let us prove by induction on $n\geq 1$
that the map $D \to B$ is onto on $p$-transitions for
$p\leq n$. 

Let $(\alpha,u_1,\dots,u_{n+1},\beta)$ be a $(n+1)$-transition of $B$,
giving rise to a map \[C_{n+1}^{ext}[\mu(u_1),\dots,\mu(u_{n+1})]
\longrightarrow B,\] which factors as a composite
\[C_{n+1}^{ext}[\mu(u_1),\dots,\mu(u_{n+1})] \to
C_{n+1}[\mu(u_1),\dots,\mu(u_{n+1})] \to B\] because $B$ is
cubical. The composite \[\de C_{n+1}[\mu(u_1),\dots,\mu(u_{n+1})]
\subset C_{n+1}[\mu(u_1),\dots,\mu(u_{n+1})] \to B\] factors
uniquely as a composite $\de C_{n+1}[\mu(u_1),\dots,\mu(u_{n+1})]
\to D \to B$ by the induction hypothesis. We obtain a
commutative diagram of cubical transition systems
\[
\xymatrix{\de C_{n+1}[\mu(u_1),\dots,\mu(u_{n+1})] \fD{} \fR{} && D \fD{} \\
&& \\
C_{n+1}[\mu(u_1),\dots,\mu(u_{n+1})] \ar@{-->}[rruu]^{k}\fR{} && B.
}
\]
The existence of the lift $k$ implies that the transition
$(\alpha,u_1,\dots,u_{n+1},\beta)$ belongs to $D$, hence the map $D
\to B$ is onto on $(n+1)$-transitions. So we obtain
$D\iso B$.

The map $C_0 \sqcup C_0 \sqcup C_1[x] \to \dd{x}$ is the composite
$C_0 \sqcup C_0 \sqcup C_1[x]\to C_0 \sqcup C_0 \sqcup \dd{x} \to
\dd{x}$ where the left-hand map is a pushout of the generating
cofibration $C_1[x] \to \dd{x}$ and where the right-hand map is a
pushout of the generating cofibration $R:\{0,1\} \to \{0\}$ twice. So
$\cell_{\cts}(\I)$ is the class of cofibrations. Therefore
$\cell_{\cts}(\I)$ is closed under retract and $\cell_{\cts}(\I) =
\cof_{\cts}(\I)$.  \epf

\bc \label{homotopy-of-cts} Let $S$ be an arbitrary set of maps in
$\cts$. There exists a unique combinatorial model category structure
on $\cts$ such that $\I$ is the set of generating cofibrations and
such that the fibrant objects are the
$\Lambda_{\cts}(V,S^{cof},\I)$-injective objects where $S^{cof}$ is a
set of cofibrant replacements of the maps of $S$. All objects are
cofibrant. The class of weak equivalences is the localizer generated
by $S$.  \ec

\bpf This corollary is a consequence of Theorem~\ref{generators} and
Theorem~\ref{build-model-cat}.  \epf

\section{Labelled symmetric precubical sets}
\label{rappel_sym_pre_set}

Let $[n] = \{0,1\}^n$ for $n \geq 0$. The unique member of the
singleton set $[0]$ is denoted by $()$. The set $[n]$ is equipped with
the partial ordering $\{0<1\}^n$.  Let $\delta_i^\alpha : [n-1] \to
[n]$ be the set map defined for $1\leq i\leq n$ and $\alpha \in
\{0,1\}$ by $\delta_i^\alpha(\epsilon_1, \dots, \epsilon_{n-1}) =
(\epsilon_1, \dots, \epsilon_{i-1}, \alpha, \epsilon_i, \dots,
\epsilon_{n-1})$.  These maps are called the \textit{face maps}.  Let
$\sigma_i:[n] \to [n]$ be the set map defined for $1\leq i\leq n-1$
and $n\geq 2$ by $\sigma_i(\epsilon_1, \dots, \epsilon_{n}) =
(\epsilon_1, \dots, \epsilon_{i-1},\epsilon_{i+1},\epsilon_{i},
\epsilon_{i+2},\dots,\epsilon_{n})$. These maps are called the
\textit{symmetry maps}. The subcategory of $\set$ generated by the
composites of face maps and symmetry maps is denoted by $\square_S$.
A presentation by generators and relations of $\square_S$ is given in
\cite[Section~6]{MR1988396}: they consist of the usual cocubical
relations, together with the Moore relations for symmetry operators
and an additional family of relations relating the face operators and
the symmetry operators. It will not be used in this paper.

\bd \cite{MR1988396} A {\rm symmetric precubical set} is a presheaf
over $\square_S$. The corresponding category is denoted by
$\square_S^{op}\set$. If $K$ is a symmetric precubical set, then let
$K_n := K([n])$ and for every set map $f:[m] \to [n]$ of $\square_S$,
denote by $f^* : K_n \to K_m$ the corresponding set map. \ed

Let $\square_S[n]:=\square_S(-,[n])$. It is called the
\textit{$n$-dimensional (symmetric) cube}. By the Yoneda lemma, one
has the natural bijection of sets
$\square_S^{op}\set(\square_S[n],K)\iso K_n$ for every precubical set
$K$. The \textit{boundary} of $\square_S[n]$ is the symmetric
precubical set denoted by $\de \square_S[n]$ defined by removing the
interior of $\square_S[n]$: $(\de \square_S[n])_k := (\square_S[n])_k$
for $k<n$ and $(\de \square_S[n])_k = \varnothing$ for $k\geq n$.  In
particular, one has $\de \square_S[0] = \varnothing$. An
\textit{$n$-dimensional} symmetric precubical set $K$ is a symmetric
precubical set such that $K_p = \varnothing$ for $p > n$ and $K_n \neq
\varnothing$. If $K$ is a symmetric precubical set, then $K_{\leq n}$
is the symmetric precubical set given by $(K_{\leq n})_p = K_p$ for $p
\leq n$ and $(K_{\leq n})_p = \varnothing$ for $p >n$.

\begin{nota} Let $f:K \to L$ be a morphism of symmetric
  precubical sets.  Let $n\geq 0$.  The set map from $K_n$ to $L_n$
  induced by $f$ will be sometimes denoted by $f_n$. \end{nota}

\begin{nota} Let $\de_i^\alpha = (\delta_i^\alpha)^*$. And let $s_i =
  (\sigma_i)^*$. \end{nota}

The precubical nerve of any topological space can be endowed with such
a structure: the $s_i$ maps are given by permuting the coordinates:
see \cite{MR1988396} again.

\bp (\cite[Proposition~A.4]{symcub}) \label{LABEL} The following data
define a symmetric precubical set denoted by $!^S\Sigma$:
\begin{itemize}
\item $(!^S\Sigma)_0=\{()\}$ (the empty word)
\item for $n\geq 1$, $(!^S\Sigma)_n=\Sigma^n$
\item $\de_i^0(a_1,\dots,a_n) = \de_i^1(a_1,\dots,a_n) =
  (a_1,\dots,\widehat{a_i},\dots,a_n)$ where the notation
  $\widehat{a_i}$ means that $a_i$ is removed
\item $s_i(a_1,\dots,a_n) =
  (a_1,\dots,a_{i-1},a_{i+1},a_i,a_{i+2},\dots,a_n)$ for $1\leq i\leq n$.
\end{itemize}
The map $!^S:\set \to \square_S^{op}\set$ yields a well-defined
functor from the category of sets to the category of symmetric
precubical sets.  \ep

\bd\label{def_lsps} A {\rm labelled symmetric precubical set (over
  $\Sigma$)} is an object of the comma category $\square_S^{op}\set
\ddownarrow !^S\Sigma$. That is, an object is a map of symmetric
precubical sets $\ell:K \to !^S\Sigma$ and a morphism is a
commutative diagram \[ \xymatrix{ K \ar@{->}[rr]\ar@{->}[rd]&& L
  \ar@{->}[ld]\\ & !^S\Sigma.&}
\]
The map $\ell$ is called the {\rm labelling map}.  The symmetric
precubical set $K$ is sometimes called the {\rm underlying symmetric
  precubical set} of the labelled symmetric precubical set. A labelled
symmetric precubical set $K \to !^S\Sigma$ will be denoted by
$\ls{K}$. And the set of $n$-cubes $K_n$ will be also denoted
by $\ls{K}_n$. \ed

The link between labelled symmetric precubical sets and process
algebra is detailed in \cite{ccsprecub} and in the appendix of
\cite{symcub}.

\begin{nota} Let $n\geq 1$. Let $a_1,\dots,a_n$ be labels of
  $\Sigma$. Let us denote by $\square_S[a_1,\dots,a_n] : \square_S[n]
  \to !^S\Sigma$ the labelled symmetric precubical set
  corresponding by the Yoneda lemma to the $n$-cube $(a_1,\dots,a_n)$.
  And let us denote by $\de\square_S[a_1,\dots,a_n] : \de\square_S[n]
  \to !^S\Sigma$ the labelled symmetric precubical set defined
  as the composite \[\xymatrix{\de\square_S[a_1,\dots,a_n] :
    \de\square_S[n] \subset \square_S[n]
    \fR{\square_S[a_1,\dots,a_n]}&& !^S\Sigma.}\] Every set can be
  identified with a sum of $0$-cubes $C_0[]$ (also denoted by $C_0$).
\end{nota}

Since colimits are calculated objectwise for presheaves, the $n$-cubes
are finitely accessible. Since the set of cubes is a dense (and hence
strong) generator, the category of labelled symmetric precubical sets
is locally finitely presentable by \cite[Theorem~1.20 and
Proposition~1.57]{MR95j:18001}. When the set of labels $\Sigma$ is the
singleton $\{\tau\}$, the category $\square_S^{op}\set \ddownarrow
!^S\{\tau\}$ is isomorphic to the category of (unlabelled) symmetric
precubical sets because $!^S\{\tau\}$ is the terminal symmetric
precubical set.

\section{The homotopy theory of labelled symmetric precubical sets}
\label{homotopy-cube}

This section is devoted to the construction of a model category
structure on the category $\square_S^{op}\set \ddownarrow !^S\Sigma$
of labelled symmetric precubical sets. Note that if $\Sigma$ is a
singleton, then the category is isomorphic to the category of
unlabelled symmetric precubical sets, and what follows applies as
well.

\bd The {\rm interval object} of $\square_S^{op}\set \ddownarrow
!^S\Sigma$ is the labelled symmetric precubical set $\ls{!^S(\Sigma\p
  \{0,1\})}$ induced by the projection map $\Sigma\p \{0,1\} \to
\Sigma$.  Let \[\cyl\ls{K} = \ls{!^S(\Sigma\p\{0,1\})}\p \ls{K}.\] \ed

Note that $\id_{!^S\Sigma} :\ls{!^S\Sigma}$ is the terminal object
$\mathbf{1}$ of $\square_S^{op}\set \ddownarrow !^S\Sigma$.  For $k
\in \{0,1\}$, denote by $\gamma^k : \ls{!^S\Sigma} \to
\ls{!^S(\Sigma\p \{0,1\})}$ the map of labelled symmetric precubical
sets induced by the composite set map $\Sigma \iso \Sigma\p \{k\}
\subset \Sigma\p \{0,1\}$. And denote by $\sigma: \ls{!^S(\Sigma\p
  \{0,1\})} \to \ls{!^S\Sigma}$ the canonical map, also induced by the
projection map $\Sigma \p \{0,1\} \to \Sigma$. Let $\gamma = \gamma^0
\sqcup \gamma^1 : \ls{!^S\Sigma} \sqcup \ls{!^S\Sigma} \to
\ls{!^S(\Sigma\p \{0,1\})}$.

If $\ls{K}$ and $\ls{L}$ are two labelled symmetric precubical sets,
then their binary product in $\square_S^{op}\set \ddownarrow
!^S\Sigma$ is the labelled symmetric precubical set $\ls{K
  \p_{!^S\Sigma} L}$. 

\bp \label{existence-right-adjoint} The interval object
$\ls{!^S(\Sigma\p \{0,1\})}$ is exponential.  \ep

\bpf Let $\ls{K}$ be a labelled symmetric precubical set.  Then
$(\cyl\ls{K})_n = K_n\p_{\Sigma^n} (\Sigma\p \{0,1\})^n \iso K_n \p
\{0,1\}^n$, i.e. $\ls{!^S(\Sigma\p\{0,1\})} \p \ls{K} \iso \ls{(K_* \p
  \{0,1\}^*)}$ with an obvious definition of the face and symmetry
maps. So the associated cylinder functor $\ls{!^S(\Sigma\p\{0,1\})} \p
-$ is colimit-preserving because the category of sets is
cartesian-closed and because colimits are calculated objectwise in the
category $\square_S^{op}\set \ddownarrow !^S\Sigma$. Since
$\square_S^{op}\set \ddownarrow !^S\Sigma$ is well-copowered by
\cite[Theorem~1.58]{MR95j:18001}, the cylinder is a left adjoint by
the dual of the Special Adjoint Functor Theorem
\cite{MR1712872}. Hence the interval object is exponential. \epf

\bd A map of labelled symmetric precubical sets
$f:\ls{K}\longrightarrow \ls{L}$ is a {\rm cofibration} if for every
$n\geq 1$, the set map $K_n \longrightarrow L_n$ is one-to-one. \ed

\bp The class of cofibrations is generated by the set 
\begin{multline*}\I = \{\de\square_S[a_1,\dots,a_n]
  \subset \square_S[a_1,\dots,a_n] \mid n\geq 1\hbox{ and
  }a_1,\dots,a_n \in \Sigma\} \\ \cup \{C:\varnothing \to
  \{0\},R:\{0,1\} \to \{0\}\},\end{multline*} i.e. the class of
cofibrations is exactly $\cof_{\square_S^{op}\set \ddownarrow
  !^S\Sigma}(\I)$. Moreover, one has \[\cell_{\square_S^{op}\set
  \ddownarrow !^S\Sigma}(\I) = \cof_{\square_S^{op}\set \ddownarrow
  !^S\Sigma}(\I).\]  \ep

\bpf This kind of proof is standard. Let $f: \ls{K} \to \ls{L}$ be a
cofibration of labelled symmetric precubical sets. Let $\I_{0} =
\{C:\varnothing \to \{0\},R:\{0,1\} \to \{0\}\}$, and for $n\geq 1$,
let $\I_n = \{\de\square_S[a_1,\dots,a_n] \subset
\square_S[a_1,\dots,a_n]\mid a_1,\dots,a_n \in \Sigma\}$. Let
$f=f_0$. Factor $f_0$ as a composite $\ls{K} \to \ls{K^1}
\stackrel{f^1}\to \ls{L}$ where the left-hand map belongs to
$\cell_{\square_S^{op}\set \ddownarrow !^S\Sigma}(\I_{0})$ and where
the right-hand map belongs to $\inj_{\square_S^{op}\set \ddownarrow
  !^S\Sigma}(\I_{0})$. Then $f^1$ is bijective on $0$-cubes and by
hypothesis is one-to-one on $n$-cubes with $n\geq 1$. Let us suppose
$f^n:\ls{K^n}\to \ls{L}$ constructed for $n\geq 1$ and let us suppose
that it is bijective on $k$-cubes for $k<n$ and one-to-one on
$k$-cubes for $k\geq n$. Consider the pushout diagram of labelled
symmetric precubical sets
\[
\xymatrix
{
\bigsqcup \de\square_S[a_1,\dots,a_n] \fR{} \fD{} && \ls{K^n} \fD{} \\
&& \\
\bigsqcup \square_S[a_1,\dots,a_n] \fR{} && \cocartesien\ls{K^{n+1}}
}
\]
where the sum is over all commutative squares of the form 
\[
\xymatrix
{
\de\square_S[a_1,\dots,a_n] \fR{} \fD{} && \ls{K^n} \fD{} \\
&& \\
\square_S[a_1,\dots,a_n] \fR{} && \ls{L}.
}
\]
Then the map $f^{n+1}:\ls{K^{n+1}} \to \ls{L}$ is bijective on
$k$-cubes for $k<n+1$ and one-to-one on $k$-cubes for $k\geq
n+1$. Hence $f=\liminj f_n$.
\epf

\begin{rem} The sets of generating cofibrations of $\square_S^{op}\set
  \ddownarrow !^S\Sigma$ and $\cts$ are both denoted by $\I$. The
  context will always enable the reader to avoid any confusion.
\end{rem}

\bp \label{wfs_OK} The codiagonal $\ls{!^S\Sigma} \sqcup
\ls{!^S\Sigma} \longrightarrow \ls{!^S\Sigma}$ factors as a composite
\[\ls{!^S\Sigma} \sqcup \ls{!^S\Sigma} \longrightarrow \ls{!^S(\Sigma\p
  \{0,1\})}\longrightarrow \ls{!^S\Sigma}\] such that the left-hand
map induces a cofibration $\ls{K} \sqcup \ls{K} \to \cyl\ls{K}$ for
any labelled symmetric precubical set $\ls{K}$.  In other terms, the
interval object $\ls{!^S(\Sigma\p \{0,1\})}$ is good.  \ep

\bpf The left-hand map is induced by the two inclusions $\Sigma \iso
\Sigma \p \{\epsilon \}\subset \Sigma\p \{0,1\}$ with $\epsilon =
0,1$.  The right-hand map is induced by the projection $\Sigma\p
\{0,1\} \longrightarrow \Sigma$. For $n\geq 1$, the left-hand map
induces on the sets of $n$-cubes the one-to-one set map $(\Sigma \p
\{0\})^n \sqcup (\Sigma \p \{1\})^n \subset (\Sigma \p \{0,1\})^n$. So
for any labelled symmetric precubical set $\ls{K}$, and any $n\geq 1$,
the map $\ls{K} \sqcup \ls{K} \to \cyl\ls{K}$ induces on the sets
of $n$-cubes the one-to-one set map \[(K_n \p \{0\}^n) \sqcup (K_n\p
\{1\}^n) \subset (K_n \p \{0,1\})^n.\] Note that the set map
$(!^S\Sigma \sqcup !^S\Sigma)_0 \longrightarrow (!^S(\Sigma\p
\{0,1\}))_0$ is not one-to-one because it is isomorphic to the set map
$R: \{0,1\} \to \{0\}$.  \epf

The interval object $\ls{!^S(\Sigma\p \{0,1\})}$ is not very good. It
is easy to prove that the right-hand map satisfies the right lifting
property with respect to all generating cofibrations except the
cofibrations $\de\square_S[x,y]\to \square_S[x,y]$ for $x,y\in
\Sigma$. Indeed, in the commutative square of solid arrows of labelled
symmetric precubical sets
\[
\xymatrix{
 \de\square_S[x,y]  \fR{g} \fD{\subset} && !^S(\Sigma\p \{0,1\}) \fD{}\\
  && \\
  \square_S[x,y] \fR{} \ar@{-->}[rruu]^-{k} && !^S\Sigma}
\] 
the lift $k$ exists if and only if two opposite faces of
$\de\square_S[x,y]$ are labelled by $g$ in $!^S(\Sigma\p \{0,1\})$ by
the same element of $\Sigma\p \{0,1\}$.

\bp \label{cartesian} For every cofibration $f : \ls{K} \to \ls{L}$ of
labelled symmetric precubical sets, the maps $f\star \gamma$ and
$f\star \gamma^\epsilon$ for $\epsilon=0,1$ are cofibrations as
well. In other terms, the interval object $\ls{!^S(\Sigma\p\{0,1\})}$
is cartesian. \ep

\bpf The map $f\star \gamma : (\ls{L} \sqcup \ls{L})
\sqcup_{\ls{K}\sqcup \ls{K}} \cyl\ls{K} \to \cyl\ls{L}$ is a
cofibration because for $n \geq 1$, the set map
\[(f\star \gamma)_n : (L_n \sqcup L_n) \sqcup_{K_n \sqcup K_n} (K_n \p
\{0,1\}^n) \to L_n \p \{0,1\}^n\] is one-to-one. Indeed, it consists
of the inclusions $K_n\subset K_n \p \{0\}^n \subset L_n \p
\{0\}^n\subset L_n \p \{0,1\}^n$ and $K_n\subset K_n \p \{1\}^n
\subset L_n \p \{1\}^n\subset L_n \p \{0,1\}^n$. The map $f\star
\gamma^\epsilon : \ls{L} \sqcup_{\ls{K}} \cyl\ls{K} \to \cyl\ls{L}$
with $\epsilon\in \{0,1\}$ is a cofibration because for $n \geq 1$,
the map (with $K_n$ embedded in $K_n\p \{\epsilon\}^n$ and $L_n$
embedded in $L_n\p \{\epsilon\}^n$)
\[(f\star \gamma^\epsilon)_n : L_n \sqcup_{K_n} (K_n \p \{0,1\}^n)
\longrightarrow L_n \p \{0,1\}^n\] is one-to-one. Indeed, it consists
of the inclusions $K_n \subset K_n\p \{\epsilon\}^n \subset L_n\p
\{\epsilon\}^n \subset L_n \p \{0,1\}^n$.  \epf

Hence the theorem: 

\bth \label{constr_model} There exists a unique combinatorial model
category structure on $\square_S^{op}\set \ddownarrow !^S\Sigma$ such
that the class of cofibrations is generated by $\I$ and such that the
fibrant objects are the $\Lambda_{\square_S^{op}\set \ddownarrow
  !^S\Sigma}(\ls{!^S(\Sigma\p\{0,1\})},\varnothing,\I)$-injective
objects. All objects are cofibrant.  \eth

\section{Realizing labelled precubical sets as cubical transition systems}
\label{realization}

We want now to recall the construction of the realization functor from
labelled symmetric precubical sets to cubical transition systems, as
expounded in Section~8 and Section~9 of \cite{hdts}. In the same way
as for the exposition of the homotopy theory of cubical transition
systems, this third paper of the series contains an improvement. We
also explain in this section how to make this functor a Quillen
functor by adding one generating cofibration to the model category of
cubical transition systems.

\begin{nota} $\cube(\square^{op}_S\set\ddownarrow !^S\Sigma)$ is the
  full subcategory of that of labelled symmetric precubical sets
  containing the labelled cubes $\square_S[a_1,\dots,a_n]$ with $n\geq
  0$ and $a_1,\dots,a_n\in \Sigma$. 
\end{nota}

\begin{nota} $\cube(\whdts)$ is the full subcategory of that of weak
  higher dimensional transition systems containing the labelled cubes
  $C_n[a_1,\dots,a_n]$ with $n\geq 0$ and with
  $a_1,\dots,a_n\in\Sigma$.
\end{nota}

The following theorem is new and is an improvement of
\cite[Theorem~8.5]{hdts}.

\bth \label{isocube} There exists one and only one functor \[\T :
\cube(\square^{op}_S\set\ddownarrow !^S\Sigma) \longrightarrow
\cube(\whdts)\] such that $\T(\square_S[a_1,\dots,a_n]) :=
C_n[a_1,\dots,a_n]$ for all $n\geq 0$ and all $a_1,\dots,a_n\in\Sigma$
and such that for any map of labelled symmetric precubical sets $f :
\square_S[a_1,\dots,a_m] \to \square_S[b_1,\dots,b_n]$, the map
$\{0,1\}^m \to \{0,1\}^n$ between the sets of states induced by
$\T(f)$ is the map induced by $f$ between the sets of
$0$-cubes. Moreover this functor yields an isomorphism of
categories \[\cube(\square^{op}_S\set\ddownarrow !^S\Sigma) \iso
\cube(\whdts).\] \eth

\bpf[Sketch of proof] Let $m,n\geq 0$ and $a_1,\dots,a_m,b_1,\dots,b_n
\in \Sigma$. A map of labelled symmetric precubical sets $f :
\square_S[a_1,\dots,a_m] \to \square_S[b_1,\dots,b_n]$ gives rise to a
set map $f_0 : \{0,1\}^m \to \{0,1\}^n$ from the set of states of
$C_m[a_1,\dots,a_m]$ to the set of states of $C_n[b_1,\dots,b_n]$
which belongs to $\square_S([m],[n]) =
\square_S^{op}\set(\square_S[m],\square_S[n])$. By
\cite[Lemma~8.1]{hdts}, there exists a unique set map $\widehat{f} :
\{1,\dots,n\} \to \{1,\dots,m\} \cup \{-\infty,+\infty\}$ such that
$f_0(\epsilon_1,\dots,\epsilon_m) =
(\epsilon_{\widehat{f}(1)},\dots,\epsilon_{\widehat{f}(n)})$ for every
$(\epsilon_1,\dots,\epsilon_m)\in [m]$ with the conventions
$\epsilon_{-\infty} = 0$ and $\epsilon_{+\infty} = 1$. Moreover, the
restriction $\overline{f} : \widehat{f}^{-1}(\{1,\dots,m\}) \to
\{1,\dots,m\}$ is a bijection. Since $f : \square_S[a_1,\dots,a_m] \to
\square_S[b_1,\dots,b_n]$ is compatible with the labelling, one
necessarily has $a_i = b_{\overline{f}^{-1}(i)}$ for every $i\in
\{1,\dots,m\}$. One deduces a set map $\widetilde{f} :
\{(a_1,1),\dots,(a_m,m)\} \to \{(b_1,1),\dots,(b_n,n)\}$ from the set
of actions of $C_m[a_1,\dots,a_m]$ to the set of actions of
$C_n[b_1,\dots,b_n]$ by setting $\widetilde{f}(a_i,i) =
(b_{\overline{f}^{-1}(i)},\overline{f}^{-1}(i)) =
(a_i,\overline{f}^{-1}(i))$. By \cite[Lemma~8.2]{hdts}, if $f :
\square_S[a_1,\dots,a_m] \to \square_S[b_1,\dots,b_n]$ and $g :
\square_S[b_1,\dots,b_n] \to \square_S[c_1,\dots,c_p]$ are two maps of
labelled symmetric precubical sets, then one has $\widehat{g\circ f} =
\widehat{f} \circ \widehat{g}$. And by \cite[Lemma~8.3]{hdts}, the two
set maps $f_0$ and $\widetilde{f}$ above defined by starting from a
map of labelled symmetric precubical sets $f :
\square_S[a_1,\dots,a_m] \to \square_S[b_1,\dots,b_n]$ yield a map of
weak higher dimensional transition systems $\T(f) : C_m[a_1,\dots,a_m]
\to C_n[b_1,\dots,b_n]$. Hence the proof is complete with
\cite[Proposition~8.4]{hdts} and \cite[Theorem~8.5]{hdts}.  \epf

\bth \cite[Theorem~9.2]{hdts} \label{colim1} There exists a unique
colimit-preserving functor
\[\T: {\square_S}^{op}\set \ddownarrow !^S\Sigma \to \whdts\]
extending the functor $\T$ previously constructed on the full
subcategory of labelled cubes. Moreover, this functor is a left
adjoint. \eth

One has $\T(\square_S^{op}\set \ddownarrow !^S\Sigma) \subset \cts$
since every colimit of cubes is cubical by
\cite[Theorem~3.11]{cubicalhdts} (see also
\cite[Proposition~9.8]{hdts}). But we have the surprising negative result:

\bp \label{not-left-Quillen} The restriction $\T:\square_S^{op}\set
\ddownarrow !^S\Sigma \to \cts$ is not a left Quillen functor. \ep

\bpf Consider the cofibration $\de \square_S[x,y] \subset
\square_S[x,y]$ where $x$ and $y$ are two elements of $\Sigma$.  The
map of cubical transition systems $\T(\de \square_S[x,y] \subset
\square_S[x,y])$ is not a cofibration of cubical transition systems
because the set of actions of $\T(\de \square_S[x,y])$ is the set with
four elements $\{x_1,x_2,y_1,y_2\}$ with $\mu(x_1) = \mu(x_2) = x$ and
$\mu(y_1) = \mu(y_2) = y$ and the set of actions of
$\T(\square_S[x,y])$ is the set with two elements $\{x,y\}$. \epf

To overcome the problem, we are going to add the maps $\T(\de
\square_S[x,y] \subset \square_S[x,y])$ with $x,y$ running over
$\Sigma$ to the set of generating cofibrations of
$\cts$. Surprisingly, the Olschok construction can be used again with
the same interval object.

\begin{nota}
Let $\I^+ = \I \cup \{\T(\de
\square_S[x,y] \subset \square_S[x,y])\mid x,y\in\Sigma\}$. 
\end{nota}

\bp \label{p2} For any labelled symmetric precubical set $\ls{K}$, one
has \[\T(\ls{!^S(\Sigma\p\{0,1\})} \p \ls{K}) = V \p \T(\ls{K}).\] \ep

\bpf One can suppose without loss of generality that $K =
\square_S[x_1,\dots,x_n]$ for $n\geq 0$ and $x_1,\dots,x_n\in \Sigma$
because all involved functors are colimit-preserving. For $n=0$, the
two members of the equality are isomorphic to the $1$-state cubical
transition system $\{()\}$. Now suppose that $n\geq 1$. The sets of
states of the two members of the equality are equal to
$\{0,1\}^n$. Since the two weak HDTS are cubical, all actions are
used. So it suffices to check that the two members of the equality
have the same set of transitions by using the fact that a cubical
transition system is the union of its subcubes by
Definition~\ref{subcube}.

By the Yoneda lemma, the $p$-cubes of $\ls{K}$ for $p \geq 1$ are in
bijection with the maps of labelled symmetric precubical sets
$\square_S[x_{\phi(1)},\dots,x_{\phi(p)}] \to \ls{K}$ where
$\phi:\{1,\dots,p\} \to \{1,\dots,n\}$ is a one-to-one map. The image
of the $p$-transition
\[(0_p,(x_{\phi(1)},\phi(1)),\dots,(x_{\phi(p)},\phi(p)),1_p)\] of
$\T(\square_S[x_{\phi(1)},\dots,x_{\phi(p)}]) =
C_p[x_{\phi(1)},\dots,x_{\phi(p)}]$ is a $p$-transition of the
form~\footnote{Remember that a $p$-cube
  $C_p[u_1,\dots,u_p]$ has the set of actions
  $\{(u_1,1),\dots,(u_p,p)\}$.}
\[((\alpha_1,\dots,\alpha_n),(x_{\phi(1)},\phi(1)),\dots,(x_{\phi(p)},\phi(p)),(\beta_1,\dots,\beta_n))\]
where $\alpha_i \leq \beta_i$ for all $i\in \{1,\dots,n\}$ with
equality if and only if $i$ does not belong to the image of
$\phi$. 

Similarly, using the calculation in the proof of
Proposition~\ref{existence-right-adjoint}, the $p$-cubes of
$\ls{!^S(\Sigma\p\{0,1\})} \p \ls{K}$ for $p \geq 1$ are in bijection
with the maps of labelled symmetric precubical sets
$\square_S[(x_{\phi(1)},\epsilon_1),\dots,(x_{\phi(p)},\epsilon_p)]
\to \ls{!^S(\Sigma\p\{0,1\})} \p \ls{K}$ where $\phi:\{1,\dots,p\} \to
\{1,\dots,n\}$ is a one-to-one map and where
$\epsilon_1,\dots,\epsilon_p \in \{0,1\}$. The image of the
$p$-transition
\[(0_p,(x_{\phi(1)},\phi(1),\epsilon_1),\dots,(x_{\phi(p)},\phi(p),\epsilon_p),1_p)\]
of
$\T(\square_S[(x_{\phi(1)},\epsilon_1),\dots,(x_{\phi(p)},\epsilon_p)])
= C_p[(x_{\phi(1)},\epsilon_1),\dots,(x_{\phi(p)},\epsilon_p)]$ is a
$p$-transition of the form
\[((\alpha_1,\dots,\alpha_n),(x_{\phi(1)},\phi(1),\epsilon_1),\dots,(x_{\phi(p)},\phi(p),\epsilon_p),(\beta_1,\dots,\beta_n))\]
where $\alpha_i \leq \beta_i$ for all $i\in \{1,\dots,n\}$ with
equality if and only if $i$ does not belong to the image of $\phi$.
Since $\ls{!^S(\Sigma\p\{0,1\})} \p \ls{K}$ is cubical, it is equal to
the union of its subcubes by Definition~\ref{subcube}. So all
transitions of $\ls{!^S(\Sigma\p\{0,1\})} \p \ls{K}$ are of the form
above.

For similar reasons, the set of transitions of the cubical transition system 
\[\T(\ls{K})\] 
consists of the tuples of the form 
\[((\alpha_1,\dots,\alpha_n),(x_{\phi(1)},\phi(1)),\dots,(x_{\phi(p)},\phi(p)),(\beta_1,\dots,\beta_n))\]
with $\phi:\{1,\dots,p\} \to \{1,\dots,n\}$ one-to-one and where
$\alpha_i \leq \beta_i$ for all $i\in \{1,\dots,n\}$ with equality if
and only if $i$ does not belong to the image of $\phi$. Hence the set
of transitions of $V \p \T(\ls{K})$ is equal to the one of
$\T(\ls{!^S(\Sigma\p\{0,1\})} \p \ls{K})$ by
\cite[Proposition~5.5]{cubicalhdts} (see also the description of the
cylinder in the proof of \cite[Proposition~5.8]{cubicalhdts}) and the
proof is complete.  \epf

\bp \label{p0} The functor $!^S:\set\to \square_S^{op}\set$ of
Proposition~\ref{LABEL} induces a well-defined functor from
$\set\ddownarrow \Sigma$ to $\square_S^{op}\set \ddownarrow
!^S\Sigma$. And one has the equality of functors \[\T \circ !^S =
\mathbb{S}:\set\ddownarrow \Sigma \to \cts.\] \ep

\bpf Let $\mu:L\to \Sigma$ be a set map. The two cubical transition
systems $(\T \circ !^S)(L\to \Sigma)$ and $\mathbb{S}(L\to \Sigma)$
have, for any set map $L \to \Sigma$, the same set of states $\{0\}$
and the same set of actions $L$. So, by definition of a cubical
transition system, the set of transitions of $(\T \circ !^S)(L\to
\Sigma)$ is a subset of $\{0\} \p \bigcup_{n\geq 1} L^n \p \{0\}$. The
latter set turns out to be the set of transitions of $!^S(L\to
\Sigma)$ by definition of $!^S$. Therefore the identity maps of
$\{0\}$ and of $L$ induce a well-defined map of cubical transition
systems $(\T \circ !^S)(L\to \Sigma) \to \mathbb{S}(L\to \Sigma)$ which
is bijective on states and actions, and one-to-one on transitions by
Proposition~\ref{transition-inj}. To complete the proof, we notice
that any tuple $(0,u_1,\dots,u_n,0)$ is a transition of the image by
$\mathbb{T}$ of the $n$-cube $(u_1,\dots,u_n)$ of $!^S(L\to \Sigma)$.
Therefore the map $(\T \circ !^S)(L\to \Sigma) \to \mathbb{S}(L\to
\Sigma)$ is surjective on transitions.  \epf

\bp \label{p1} One has $\T(\ls{!^S(\Sigma\p\{0,1\})}) = V$. Moreover,
for $k=0,1$, the image by $\T$ of the map of labelled symmetric
precubical sets $\gamma^k:\mathbf{1} \to \ls{!^S(\Sigma\p\{0,1\})}$
induced by $\Sigma \iso \Sigma\p\{k\} \subset \Sigma\p \{0,1\}$ is the
map of cubical transition systems $\mathbf{1} \to V$ induced by
$\Sigma \iso \Sigma\p\{k\} \subset \Sigma\p \{0,1\}$.  \ep

\bpf The equality $\T(\ls{!^S(\Sigma\p\{0,1\})}) = V$ comes from
Proposition~\ref{p0} applied to the projection map $\Sigma \p\{0,1\}
\to \Sigma$. The last part about $\gamma^k$ is a corollary of
Proposition~\ref{p0} applied to the set map $\Sigma \iso \Sigma\p\{k\}
\subset \Sigma\p \{0,1\}$. \epf

We can now introduce the ``augmented'' model category structure of
cubical transition systems.

\bth \label{homotopy-of-cts-plus} There exists a unique combinatorial
model category structure on $\cts$ such that the set of maps $\I^+$ is
the set of generating cofibrations and such that the fibrant objects
are the $\Lambda_{\cts}(V,\varnothing,\I^+)$-injective objects. All
objects are cofibrant. This model category will be denoted by
$\cts^+$. \eth

\bpf The interval object $V$ is still good~\footnote{The interval
  object $V$ is not very good with respect to the maps of $\I^+$
  because there is no lift $k$ in the following diagram with $x,y\in \Sigma$ 
if the set of actions of $\T(\de\square_S[x,y])$ is taken to be 
$\{(x,0),(x,1),(y,0),(y,1)\}$.
\[
\xymatrix
{
\T(\de\square_S[x,y]) \fr{} \fd{} & \T\ls{!^S(\Sigma\p\{0,1\})} \\
\T(\square_S[x,y]) \ar@{-->}[ru]^{k}&
}
\]
The lift $k$
would exist if and only if two opposite faces of the empty square
$\T(\de\square_S[x,y])$ were labelled by the same
action.\label{not-good}} with respect to the maps of
$\I^+$. We already know that the interval object $V$
is exponential. To prove that it is cartesian with respect to $\I^+$,
we just have to prove that for any $x,y\in \Sigma$, $\T(f_{x,y})\star
\gamma^k$ and $\T(f_{x,y})\star \gamma$ belong to $\cof_{\cts}(\I^+)$
for $k=0,1$ and where $f_{x,y}:\de\square_S[x,y] \subset
\square_S[x,y]$ is the inclusion.  By Proposition~\ref{cartesian},
$f_{x,y}\star \gamma^k$ and $f_{x,y}\star \gamma$ are cofibrations of
$\square_S^{op}\set \ddownarrow !^S\Sigma$. By Proposition~\ref{p2}
and Proposition~\ref{p1}, $\T(f_{x,y}\star \gamma^k) = \T(f_{x,y})
\star \gamma^k$ and $\T(f_{x,y}\star \gamma) = \T(f_{x,y}) \star
\gamma$, hence the interval object $V$ is cartesian with respect to
$\I^+$.  All cubical transition systems are still cofibrant for this
new class of cofibrations. Hence the result by
Theorem~\ref{build-model-cat}.  \epf

\bd \label{csa1} A cubical transition system satisfies CSA1 (the {\rm
  First Cattani-Sassone axiom}) if for every transition
$(\alpha,u,\beta)$ and $(\alpha,u',\beta)$ such that the actions $u$
and $u'$ have the same label in $\Sigma$, one has $u = u'$. \ed

The full subcategory of cubical transition systems satisfying CSA1 is
reflective by \cite[Corollary~5.7]{hdts}. The reflection is denoted by
$\CSA_1:\cts \to \cts$.

\bp \label{ortho-fibrant} Let $T$ be a cubical transition system
satisfying CSA1. Let $I$ and $S$ be arbitrary sets of maps of
$\cts$. Then $T$ is $\Lambda_{\cts}(V,S,I)$-injective if and only if
it is $S$-orthogonal.  \ep

\bpf In the case in which $I=\I$, the set of generating cofibrations
for $\cts$, this is \cite[Proposition~7.7]{cubicalhdts}, and the same
proof works for arbitrary sets.  \epf

\bp \label{CSA1-fibrant} Every cubical transition system satisfying
CSA1 is fibrant both in $\cts$ and in $\cts^+$. \ep

\bpf The first part is \cite[Proposition~7.8]{cubicalhdts}. Every
cubical transition system satisfying CSA1 is $\varnothing$-orthogonal,
so fibrant in $\cts^+$ as well by Proposition~\ref{ortho-fibrant} and
Theorem~\ref{homotopy-of-cts-plus}. \epf

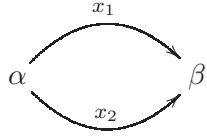
\begin{figure}
\[
\xymatrix{
\alpha \ar@/^20pt/[rr]^-{x_1} \ar@/_20pt/[rr]^-{x_2}&& \beta }
\]
\caption{Cubical transition system $\cyl(C_1[x])$ homotopy equivalent
  to the $1$-cube $C_1[x]$ with $\mu(x_1) = \mu(x_2) = x$.}
\label{Da1}
\end{figure}

\bp \label{CSA1-replacement} For every cubical transition system $X$,
the unit $X \to \CSA_1(X)$ is a weak equivalence of both $\cts$ and
$\cts^+$. \ep

\bpf The argument of the proof of \cite[Theorem~7.10]{cubicalhdts}
must be used: the unit $X \to \CSA_1(X)$ is a transfinite composition
of pushouts of maps of the form $\sigma_{C_1[x]}:\cyl(C_1[x]) \to
C_1[x]$ (the source is depicted in Figure~\ref{Da1}) with $x\in
\Sigma$. Consider such a pushout:
\[
\xymatrix
{
\cyl(C_1[x]) \fD{\sigma_{C_1[x]}} \fR{\phi} &&  X \fD{f}\\
&& \\
C_1[x] \fR{} && \cocartesien Y.
}
\] 
The map $\sigma_{C_1[x]}$ is \emph{never} a cofibration of course. The
point is that $\phi$ is either a cofibration of cubical transition
systems or it takes the two actions of $\cyl(C_1[x])$ to the same
action of $X$: in the first case, $f$ is a weak equivalence in $\cts$
and $\cts^+$ because of the left properness; in the second case, $f$
is just an isomorphism.  \epf

\bth \label{carac-w} The two model categories $\cts$ and $\cts^+$ have
the same class of weak equivalences: the class of maps of cubical
transition systems $f$ such that $\CSA_1(f)$ is an isomorphism. \eth

\bpf Two maps $f,g:X\to Y$ with $Y$ satisfying CSA1 are homotopic with
respect to the cylinder $V\p -$ if and only if they are equal by
\cite[Proposition~7.4]{cubicalhdts}. Therefore two weakly equivalent
cubical transition systems in $\cts$ or $\cts^+$ satisfying CSA1 are
isomorphic: this is \cite[Proposition~7.9]{cubicalhdts}. The proof is
complete with Proposition~\ref{CSA1-replacement}.  \epf

\bp \label{almost-T-left-Quillen} The functor $\T : \square_S^{op}\set
\ddownarrow !^S\Sigma \to \cts$ preserves weak equivalences. \ep

\bpf Let $f:\ls{K}\to \ls{L}$ be a weak equivalence of labelled
symmetric precubical sets. Let $T$ be a fibrant object of the left
determined model category structure of $\cts$.  We have to prove that
the set map $\cts(\T(\ls{L}),T)/\!\!\simeq \to
\cts(\T(\ls{K}),T)/\!\!\simeq$ induced by $f$ is a bijection of sets
where $\simeq$ means the homotopy relation induced by the cylinder of
$\cts$. By adjunction, we have to check that the set
map \[(\square_S^{op}\set \ddownarrow
!^S\Sigma)(\ls{L},\R(T))/\!\!\simeq \longrightarrow
(\square_S^{op}\set \ddownarrow !^S\Sigma)(\ls{K},\R(T))/\!\!\simeq\]
induced by $f$ is a bijection of sets where $\R$ is the right adjoint
to $\T$.  Since $f$ is a weak equivalence of labelled symmetric
precubical sets, it suffices to check that $\R(T)$ is
$\Lambda_{\square_S^{op}\set \ddownarrow
  !^S\Sigma}(\ls{!^S(\Sigma\p\{0,1\})},\varnothing,\I)$-injective.  By
Proposition~\ref{p2} and Proposition~\ref{p1}, this is equivalent to
$T$ being $\Lambda_{\cts}(V,\varnothing,\I)$-injective, which holds
because $T$ is fibrant.  \epf

\bc \label{same_with_plus} There exists a zig-zag of left Quillen
functors
\[\xymatrix{
\square_S^{op}\set \ddownarrow !^S\Sigma \fR{\T} && \cts^+ && \fL{\id_{\cts}} \cts.
}\] 
Moreover, the right-hand left Quillen functor induces a Quillen equivalence
$\cts \simeq \cts^+$.  \ec

\bpf The functor $\T : \square_S^{op}\set \ddownarrow !^S\Sigma \to
\cts^+$ is a left Quillen functor by Proposition~\ref{carac-w} and
Proposition~\ref{almost-T-left-Quillen} and by definition of the
cofibrations of $\cts^+$.  The functor $\CSA_1:\cts\to\cts$ is a
cofibrant-fibrant replacement for the two model categories $\cts$ and
$\cts^+$ by Proposition~\ref{CSA1-fibrant} and by
Proposition~\ref{CSA1-replacement}. They have the same class of weak
equivalences so the left Quillen functor $\cts\to\cts^+$ induces an
adjoint equivalence of categories between the homotopy categories.
\epf

\section{Homotopical property of the realization functor} 
\label{homotopy-realization}

We are now ready to compare labelled symmetric precubical sets and
cubical transition systems from a homotopy point of view.

\begin{nota} Let $\mathcal{S} = \{p_x:C_1[x] \sqcup C_1[x] \to
  \dd{x}\mid x\in \Sigma\}$. Let $\mathcal{S}^{cof} = \{p_x^{cof}\mid
  x\in \Sigma\}$ where $(-)^{cof}$ is a cofibrant replacement in
  $\cts$.
\end{nota}

Let $X$ be a cubical transition system. Let us factor in $\cts$ the
canonical map $X\to \mathbf{1}$ as a composite $X\to
\bl_{\mathcal{S}}(X) \to \mathbf{1}$ where the left-hand map belongs
to $\cell_{\cts}(\mathcal{S})$ and the right-hand map belongs to
$\inj_{\cts}(\mathcal{S})$. The functor $\bl_{\mathcal{S}} : \cts \to
\cts$ is studied in \cite{cubicalhdts}.  The next proposition explains
the image of the functor $\bl_{\mathcal{S}}$ on objects.

\bp \label{simpl} For a cubical transition system $X = (S,\mu:L\to
\Sigma,T=\bigcup_{n\geq 1}T_n)$, the following statements are equivalent:
\begin{enumerate}
\item The labelling map $\mu$ is one-to-one.
\item $X$ is $\mathcal{S}$-injective.
\item $X$ is $\mathcal{S}$-orthogonal.
\end{enumerate}
If one of these statements is true, then $X$ satisfies CSA1.  \ep

\bpf The equivalence $(2) \Leftrightarrow (3)$ is due to the fact that
all maps of $\mathcal{S}$ are epimorphisms. We have to prove now that
$(1) \Leftrightarrow (2)$.  Let us suppose $(2)$.  Let $x_1$ and $x_2$
be two actions of $X$ with $\mu(x_1) = \mu(x_2) = x$. Since $X$ is
injective with respect to $\underline{x_i} \to C_1[x]$ for $i=1,2$,
the two maps $\underline{x_i} \subset X$ factors as a composite
$\underline{x_i} \to C_1[x] \to X$. Hence there is a map $C_1[x]\sqcup
C_1[x] \to X$ sending one action of the source to $x_1$ and the other
one to $x_2$. By hypothesis, $X$ is $p_x$-injective. Therefore the
latter map factors as a composite $C_1[x]\sqcup C_1[x] \to \dd{x} \to
X$. So $x_1 = x_2$, and $\mu$ is one-to-one. Conversely, suppose that
$\mu$ is one-to-one. Let $C_1[x]\sqcup C_1[x] \to X$ be a map. Then
the two actions of the source are taken to the same action of $X$
because they have the same labelling. Therefore the map factors as a
composite $C_1[x]\sqcup C_1[x] \to \dd{x} \to X$.  The last assertion
is obvious. \epf

So the functor $\bl_{\mathcal{S}} : \cts \to \cts$ induces a functor
from $\cts$ to the full reflective subcategory $\mathcal{S}^\perp$ of
cubical transition systems consisting of $\mathcal{S}$-orthogonal
objects. By \cite[Theorem~8.11]{cubicalhdts}, the functor
$\bl_{\mathcal{S}}$ is left adjoint to the inclusion functor
$\iota_{\mathcal{S}} : \mathcal{S}^\perp \subset \cts$. 

\bp \label{Ls-fibrant} For every cubical transition system $X$, the
cubical transition system $\bl_{\mathcal{S}}(X)$ is
$\mathcal{S}^{cof}$-orthogonal.  \ep

\bpf We want to describe explicitly a cofibrant replacement of
$p_x:C_1[x]\sqcup C_1[x] \to \dd{x}$ in $\cts$.  The functor $V\p -$
is described in \cite[Proposition~5.5]{cubicalhdts} and in
\cite[Proposition~5.8]{cubicalhdts}. The cubical transition system $V
\p C_1[x]$ has the same state as $C_1[x]$ (one initial state $\alpha$
and one final state $\beta$), has two actions $x_1$ and $x_2$ labelled
by $x$ and two $1$-transitions $(\alpha,x_1,\beta)$ and
$(\alpha,x_2,\beta)$ (Figure~\ref{Da1}).  A cofibrant replacement of
$p_x$ can then be obtained by considering the composite map
\[p_x^{cof}: C_1[x] \sqcup C_1[x] \to (V \p C_1[x] \sqcup V\p C_1[x])
\to (V \p C_1[x] \sqcup V\p C_1[x])/(x_2=x_4=x)\] where $\{x_1,x_2\}$
is the set of actions of the left-hand copy of $V\p C_1[x]$ and where
$\{x_3,x_4\}$ is the set of actions of the right-hand copy of $V\p
C_1[x]$. To completely determine this map, we must say that it induces
a bijection on the set of states and it is the inclusion $\{x_1,x_3\}
\subset \{x_1,x_3,x\}$ on actions with $\mu(x_1) = \mu(x_2) = \mu(x_3)
= \mu(x_4) = x$. So $p_x^{cof}$ is a cofibration. The target of
$p_x^{cof}$ is depicted in Figure~\ref{Da}. One has the equalities of
cubical transition systems
\[\CSA_1\left[(V \p C_1[x] \sqcup V\p C_1[x])/(x_2=x_4=x)\right] = \dd{x} = \CSA_1(\dd{x}).\]
Therefore, using Theorem~\ref{carac-w}, the cubical transition systems
$(V \p C_1[x] \sqcup V\p C_1[x])/(x_2=x_4=x)$ and $\dd{x}$ are weakly
equivalent in $\cts$. Thus, $p_x^{cof}$ is a cofibrant replacement
of $p_x$ in $\cts$.

We now want to check that $\bl_{\mathcal{S}}(X)$ is
$p_x^{cof}$-orthogonal for all $x\in \Sigma$ to complete the
proof. Consider a commutative diagram of solid arrows:
\[
\xymatrix{
C_1[x] \sqcup C_1[x] \fD{p_x^{cof}}\fR{} && \bl_{\mathcal{S}}(X) \fD{} \\
&& \\
(V \p C_1[x] \sqcup V\p C_1[x])/(x_2=x_4=x) \ar@{-->}[rruu]^-{k} \fR{} && \mathbf{1}
}
\]
Remember by Proposition~\ref{simpl} that the labelling map of
$\bl_{\mathcal{S}}(X)$ is one-to-one. The existence and uniqueness of
the lift $k$ is then clear on states (because $p_x^{cof}$ is bijective
on states) and on actions ($k$ takes all actions of its source to
$x$). \epf

\bp \label{fibrantS1} If $X$ is a cubical transition system, then
$\bl_{\mathcal{S}}(X)$ is fibrant in $\bl_{\mathcal{S}}(\cts)$ and in
$\bl_{\mathcal{S}}(\cts^+)$. \ep

Note that the argument given in \cite[Theorem~8.11]{cubicalhdts} to
prove the fibrancy of $\bl_{\mathcal{S}}(X)$ in
$\bl_{\mathcal{S}}(\cts)$ is wrong: $\mathcal{S}$-orthogonality was
used instead of $\mathcal{S}^{cof}$-orthogonali\-ty.

\bpf A cubical transition system of the form $\bl_{\mathcal{S}}(X)$ is
$\mathcal{S}^{cof}$-orthogonal by Proposition~\ref{Ls-fibrant} and
satisfies CSA1 by Proposition~\ref{simpl}. So it is fibrant in
$\bl_{\mathcal{S}}(\cts)$ by Proposition~\ref{ortho-fibrant} and by
Corollary~\ref{homotopy-of-cts}.  The map $p_x^{cof}:C_1[x]\sqcup
C_1[x]\to B$ introduced in the proof of Proposition~\ref{Ls-fibrant}
is a cofibration of $\cts$ such that $\CSA_1(B)=\CSA_1(\dd{x})$. So it
is a cofibration of $\cts^+$ as well since there are more cofibrations
in $\cts^+$ than in $\cts$ and $B$ is weakly equivalent to $\dd{x}$ in
$\cts^+$ as well by Theorem~\ref{carac-w}. Hence $p_x^{cof}$ is a
cofibrant replacement of $p_x$ in $\cts^+$ as well and $\bl_{\mathcal{S}}(X)$ is fibrant 
in $\cts^+$ as well.  \epf

\begin{figure}
\[
\xymatrix{
0 \ar@/^20pt/[rr]^-{x_1} \ar@/_20pt/[rr]^-{x}&& 1 }
\xymatrix{
2 \ar@/^20pt/[rr]^-{x_3} \ar@/_20pt/[rr]^-{x}&& 3 }
\]
\caption{$(V \p C_1[x] \sqcup V\p C_1[x])/(x_2=x_4=x)$.}
\label{Da}
\end{figure}

\bp \label{S-replacement} For every cubical transition system $X$, the
unit $X \to \bl_{\mathcal{S}}(X)$ is a weak equivalence of both
$\bl_{\mathcal{S}}(\cts)$ and $\bl_{\mathcal{S}}(\cts^+)$. \ep

\bpf The map $X \to \bl_{\mathcal{S}}(X)$ is a transfinite composition
of pushouts of maps of $\mathcal{S}$. Then we can use the argument
\cite[Proposition~8.5]{cubicalhdts}. Let us briefly recall
it. Consider a pushout in $\cts$ of the form (with $x\in \Sigma$)
\[
\xymatrix
{
C_1[x] \sqcup C_1[x] \fD{p_x} \fR{\phi} &&  X \fD{f}\\
&& \\
\dd{x} \fR{} && \cocartesien Y.
}
\] 
The map $p_x$ is \emph{never} a cofibration of course. The point is
that $\phi$ is either a cofibration or it takes the two actions of
$C_1[x] \sqcup C_1[x]$ to the same action of $X$: in the first case,
$f$ is a weak equivalence in $\bl_{\mathcal{S}}(\cts)$ and
$\bl_{\mathcal{S}}(\cts^+)$ because of the left properness; in the
second case, $f$ is just an isomorphism.  \epf

\bth \label{carac-wS} The two model categories
$\bl_{\mathcal{S}}(\cts)$ and $\bl_{\mathcal{S}}(\cts^+)$ have the
same class of weak equivalences: the class of maps of cubical
transition systems $f$ such that $\bl_{\mathcal{S}}(f)$ is an
isomorphism. \eth

\bpf By Proposition~\ref{fibrantS1} and
Proposition~\ref{S-replacement}, the functor $\bl_{\mathcal{S}}:\cts
\rightarrow \cts$ is a cofibrant-fibrant replacement both in
$\bl_{\mathcal{S}}(\cts)$ and $\bl_{\mathcal{S}}(\cts^+)$. A map $f$
is a weak equivalence in the Bousfield localization if and only if
$\bl_{\mathcal{S}}(f)$ is a weak equivalence of $\cts$ (or of
$\cts^+$) by \cite[Theorem~3.2.13]{ref_model2}. But the source and
target of $\bl_{\mathcal{S}}(f)$ satisfy CSA1 by
Proposition~\ref{simpl}, hence the conclusion using
Theorem~\ref{carac-w}.  \epf

\bc \label{same_with_plus_and_S} The left Quillen functor
$\bl_{\mathcal{S}}(\cts) \to \bl_{\mathcal{S}}(\cts^+)$ induced by the
identity functor is a Quillen equivalence. \ec

\bpf The left Quillen functor $\bl_{\mathcal{S}}(\cts) \to
\bl_{\mathcal{S}}(\cts^+)$ induces an adjoint equivalence of
categories between the homotopy categories.  \epf

Before stating the next theorem, we need to introduce again a few
notations. 

\bd \label{def-cub} Let $X\in \whdts$. The {\rm cubification functor}
is the functor
\[\cub : \whdts \longrightarrow \whdts\] defined by \[\cub(X) =
\liminj_{C_n[x_1,\dots,x_n] \to X} C_n[x_1,\dots,x_n],\] the colimit
being taken in $\cts$ (or equivalently in $\whdts$).  \ed

We can now prove one of the main result of this paper:

\bth \label{same1} There exists a Bousfield localization of
$\square_S^{op}\set \ddownarrow !^S\Sigma$ which is Quill\-en equivalent
to the Bousfield localization of the left determined model category
structure of $\cts$ by the cubification functor.  \eth

\bpf By \cite[Corollary~8.7]{cubicalhdts}, the Bousfield localization
of the left determined model category structure of $\cts$ by the
cubification functor is $\bl_{\mathcal{S}}(\cts)$.  By
Corollary~\ref{same_with_plus_and_S} and
\cite[Proposition~3.2]{MR1870516}, it suffices to prove that the left
Quillen functor $\square_S^{op}\set \ddownarrow !^S\Sigma \to
\bl_{\mathcal{S}}(\cts^+)$ induced by $\T$ is homotopically surjective
in the sense of \cite[Definition~3.1]{MR1870516}.

The right adjoint $\R:\cts^+ \to \square_S^{op}\set \ddownarrow
!^S\Sigma$ of $\T: \square_S^{op}\set \ddownarrow !^S\Sigma \to
\cts^+$ may be defined as follows.  Let $X$ be a cubical transition
system. The set $\R(X)$ of $n$-cubes labelled by $(a_1,\dots,a_n)\in
\Sigma^n$ is the set of maps of cubical transition systems
$C_n[a_1,\dots,a_n] \to X$ (so for $n=0$, it is the set of states)
with an obvious definition of the face maps and symmetry maps. Using
the isomorphism of categories $\cube(\square^{op}_S\set\ddownarrow
!^S\Sigma) \iso \cube(\whdts)$ by Theorem~\ref{isocube}, one deduces
that $\T(\R(X)) \iso \cub(X)$. Hence by
\cite[Proposition~8.4]{cubicalhdts} and by
\cite[Proposition~8.5]{cubicalhdts}, for any cubical transition system
$X$, the counit $\T(\R(X)) \to X$ is a weak equivalence of
$\bl_{\mathcal{S}}(\cts)$, and therefore of
$\bl_{\mathcal{S}}(\cts^+)$ by
Theorem~\ref{carac-wS}~\footnote{\cite[Proposition~8.4]{cubicalhdts}
  actually proves that the counit $\cub(X) \to X$ is a
  transfinite composition of pushouts of the maps $p_x:C_1[x]\sqcup
  C_1[x] \to \dd{x}$ for $x\in \Sigma$; so $\cub(X) \to X$ is
  a weak equivalence in the Bousfield localizations by the same
  argument as for proving Proposition~\ref{S-replacement}.}. Since all
cubical transition systems are cofibrant, we see that the last
assertion is nothing else but the definition of homotopically
surjective.  \epf

\section{The higher dimensional automata paradigm}
\label{paradigm}

If $K = \square_S[x_1,\dots,x_n] \sqcup_{\de \square_S[x_1,\dots,x_n]}
\square_S[x_1,\dots,x_n]$ with $n\geq 2$ and $x_1,\dots,x_n\in
\Sigma$, then $\R(\bl_{\mathcal{S}}(\T(K))) =
\square_n[x_1,\dots,x_n]$ by Proposition~\ref{simpl} and
\cite[Proposition~9.3]{hdts}. This suggests to recall now the HDA
paradigm for a slight improvement of Theorem~\ref{same1}.

\bd \cite[Definition~7.1]{hdts} A labelled symmetric precubical set
\[\ls{K}\] satisfies {\rm the paradigm of higher dimensional automata (HDA
  paradigm)} if for every $p\geq 2$, every commutative square of solid
arrows (called a {\rm labelled $p$-shell} or {\rm labelled
  $p$-dimensional shell})
\[
\xymatrix{
  \de\square_S[p] \fR{}\fD{} && K \fD{}\\
  &&\\
  \square_S[p] \fR{} \ar@{-->}[rruu]^-{k} && !^S\Sigma}\] admits at
most one lift $k$ (i.e. a map $k$ making the two triangles
commutative).  \ed

The interest of the HDA paradigm in computer science is that it is
satisfied by all real examples (see for example
\cite[Theorem~5.2]{ccsprecub} and \cite[Corollary~5.3]{ccsprecub}).  A
full $n$-cube with $n\geq 2$ models the concurrent execution of $n$
actions. An empty $n$-cube with $n\geq 2$ models the concurrent
execution of $n-1$ actions maximum among a set of $n$ actions
\cite{ccsprecub}. It is impossible to have two $n$-cubes (for $n\geq
2$) with the same boundary.  Either it is possible for the $n$ actions
to run concurrently (full), or there is an obstruction (empty).

Note that the HDA paradigm is automatically satisfied by higher
dimensional transition systems because for $n\geq 2$, there is the
isomorphism
\[C_n[x_1,\dots,x_n] \sqcup_{\de C_n[x_1,\dots,x_n]}
C_n[x_1,\dots,x_n] \iso C_n[x_1,\dots,x_n]\] for all $x_1,\dots,x_n
\in \Sigma$ by \cite[Proposition~9.3]{hdts}.

By \cite[Proposition~7.3]{hdts}, the HDA paradigm is equivalent to
being orthogonal to the set of maps \[\left\{\square_S[a_1,\dots,a_p]
\sqcup_{\de\square_S[a_1,\dots,a_p]}
\square_S[a_1,\dots,a_p]\to\square_S[a_1,\dots,a_p]\right\}\]
for $p\geq 2$ and $a_1, \dots, a_p\in \Sigma$. So (see
\cite[Corollary~7.4]{hdts}), the full subcategory, denoted by
$\hda^\Sigma$, of $\square_S^{op}\set \ddownarrow !^S\Sigma$
containing the objects satisfying the HDA paradigm is a full
reflective locally presentable category of the category
$\square_S^{op}\set \ddownarrow !^S\Sigma$ of labelled symmetric
precubical sets.  In other terms, the inclusion functor $i_{\Sigma}:
\hda^\Sigma \subset \square_S^{op}\set \ddownarrow !^S\Sigma$ has a
left adjoint $\sh_{\Sigma}: \square_S^{op}\set \ddownarrow !^S\Sigma
\to \hda^\Sigma$.

In fact the category $\hda^\Sigma$ is locally finitely presentable;
indeed, the labelled $n$-cubes for $n\geq 0$ are in $\hda^\Sigma$ by
\cite[Proposition~7.2]{hdts}, and one can prove that they form a dense
set of generators.

\begin{nota} When $\Sigma$ is the singleton $\{\tau\}$, let
  $i:=i_\Sigma$ and $\sh:= \sh_\Sigma$. \end{nota}

One has \[i_\Sigma\ls{K} \iso \ls{i(K)} = \ls{K}\] and \[\sh_\Sigma\ls{K}
\iso (\sh(K)\to \sh(!^S \Sigma) \iso !^S \Sigma)\] because the symmetric
precubical set $!^S \Sigma$ belongs to $\hda$.

We want to restrict the homotopy theory of labelled symmetric
precubical sets to the full reflective subcategory $\hda^\Sigma$. So
we must explain how to restrict the Olschok construction to a full
reflective subcategory, at least within our particular setting.

\bth (\cite[Lemma~5.2]{MO} with some additional
remarks) \label{restriction-model-structure} (Restriction to a full
reflective subcategory) Let $\K$, $I$, $S$ and $V$ be as in
Theorem~\ref{build-model-cat}. Let $\mathcal{A}$ be a full reflective
locally presentable subcategory with $V\in \mathcal{A}$ and $I \subset
\Mor(\mathcal{A})$.  Let $R:\K \to \mathcal{A}$ be the
reflection. Suppose that $S \subset \cof_{\mathcal{A}}(I)$ and that
for every object $X\in \mathcal{A}$, $X^V \in \mathcal{A}$ where
$(-)^V$ is the right adjoint in $\K$ to the cylinder functor
$\cyl(-)=V\p -$. Then there exists a unique combinatorial model
category structure such that the class of cofibrations is generated by
$I$ and such that an object of $\mathcal{A}$ is fibrant if and only if
it is $\Lambda_\K(V,S,I)$-injective. In particular, the fibrant
objects of $\mathcal{A}$ are the fibrant objects of $\K$ belonging to
$\mathcal{A}$. The reflection $R:\K \to \mathcal{A}$ is a
homotopically surjective left Quillen adjoint. All objects are
cofibrant.  \eth

\bpf The inclusion functor $\mathcal{A} \subset \K$ is a right
adjoint, therefore it preserves binary products. And $V\in
\mathcal{A}$ by hypothesis. Hence $V\p -:\K \to \K$ restricts to an
endofunctor of $\mathcal{A}$.

Note that $I\subset \mathcal{A}$ implies that $I = R(I)$. One has
$\inj_{\mathcal{A}}(I) = \Mor(\mathcal{A}) \cap \inj_\K(I)$ because
$\mathcal{A}$ is a full subcategory. So for all $f\in
\inj_{\mathcal{A}}(I)$, $(\varnothing \to A) \square f$ for every
object $A$ of $\mathcal{A}$ because $(\cof_\K(I),\inj_\K(I))$ is a
cofibrant small weak factorization system, hence
$(\cof_{\mathcal{A}}(I),\inj_{\mathcal{A}}(I))$ is a cofibrant small
weak factorization system of $\mathcal{A}$.

Since $R$ is colimit-preserving, $R(\cell_\K(I)) \subset
\cell_\mathcal{A}(I)$. Every map of $\cof_\K(I)$ is a retract of a map
of $\cell_\K(I)$, therefore $R(\cof_\K(I)) \subset
\cof_{\mathcal{A}}(I)$. Let $f\in \Mor(\mathcal{A}) \cap \cof_\K(I)$.
Then $R(f) = f \in \cof_{\mathcal{A}}(I)$. One obtains the inclusion
$\Mor(\mathcal{A}) \cap \cof_\K(I) \subset \cof_{\mathcal{A}}(I)$.
Let $A\in \mathcal{A}$. Then $A\sqcup A \to V\p A$ is a map of
$\mathcal{A}$ belonging to $\cof_\K(I)$ because $V$ is good in
$\K$. So $V$ is good in $\mathcal{A}$ as well.

Let $A,B \in \mathcal{A}$.  Then $\mathcal{A}(V\p A,B) \iso \K(V\p
A,B) \iso \K(A,B^V) \iso \mathcal{A}(A,B^V)$ because $B^V \in
\mathcal{A}$ and because $\mathcal{A}$ is a full subcategory of
$\K$. So the interval object $V$ is exponential in
$\mathcal{A}$. Since $R$ is colimit-preserving, one has $I \star
\gamma \in R(\cof_\K(I)) \subset \cof_\mathcal{A}(I)$ and $I \star
\gamma^k \in R(\cof_\K(I)) \subset \cof_\mathcal{A}(I)$, hence $V$ is
cartesian in $\mathcal{A}$ as well.

Let $T$ be an object of $\mathcal{A}$. Then $T$ is
$\Lambda_\K(V,S,I)$-injective, if and only if it is
$\Lambda_\mathcal{A}(V,S,I)$-injective because $R(\Lambda_\K(V,S,I)) =
\Lambda_\mathcal{A}(V,S,I)$. So the fibrant objects of $\mathcal{A}$
are the fibrant objects of $\K$ belonging to $\mathcal{A}$. Hence the
existence of the model category structure on $\mathcal{A}$ by
Theorem~\ref{build-model-cat}.

Since $R(\cof_\K(I)) \subset \cof_\mathcal{A}(I)$, the reflection
$R$ takes cofibrations to cofibrations. Let $f$ be a weak equivalence
of $\K$. Let $T$ be a fibrant object of $\mathcal{A}$. Then
$\K(f,T)/\!\!\simeq$ is a bijection. By adjunction, one gets that
$\mathcal{A}(R(f),T)/\!\!\simeq$ is a bijection as well, i.e. that
$R(f)$ is a weak equivalence of $\mathcal{A}$. Therefore the
reflection $R : \K \to \mathcal{A}$ is a left Quillen adjoint. It is
homotopically surjective because for any fibrant object $A\in
\mathcal{A}$, $A$ is cofibrant in $\K$ and there is the isomorphism
$R(A)\iso A$.  \epf

Denote by $\ls{K} \mapsto \ls{K}^{\ls{!^S(\Sigma\p\{0,1\})}}$ the
right adjoint to the cylinder functor of $\square_S^{op}\set
\ddownarrow !^S\Sigma$ which exists by
Proposition~\ref{existence-right-adjoint}.

\begin{figure}
\[
\xymatrix{
  \de\square_S[a_1,\dots,a_p] \fR{f}\fD{} && \ls{K}^{\ls{!^S(\Sigma\p\{0,1\})}}\\
  &&\\
  \square_S[a_1,\dots,a_p] \ar@{-->}[rruu]^-{\overline{f}} && }\]
\caption{$\ls{K}^{\ls{!^S(\Sigma\p\{0,1\})}}$ must satisfy the HDA paradigm.}
\label{fig1}
\end{figure}

\begin{figure}
\[
\xymatrix{
  \cyl(\de\square_S[a_1,\dots,a_p])  \fR{f}\fD{} && K\\
  &&\\
  \cyl(\square_S[a_1,\dots,a_p]). \ar@{-->}[rruu]^-{\overline{f}} && }\] 
\caption{HDA paradigm for $\ls{K}^{\ls{!^S(\Sigma\p\{0,1\})}}$ and adjunction.}
\label{fig2}
\end{figure}

\bp For every labelled symmetric precubical set $\ls{K}$, if $\ls{K}$
satisfies the HDA paradigm, then $\ls{K}^{\ls{!^S(\Sigma\p\{0,1\})}}$
satisfies the HDA paradigm as well. \ep

\bpf We have to prove that there exists at most one lift $k$ for every
commutative diagram of the form of Figure~\ref{fig1} with $p\geq 2$
and $a_1,\dots,a_p \in \Sigma$. By adjunction, that means that we have
to prove that there exists at most one list $\overline{f}$ for every
commutative diagram of the form of Figure~\ref{fig2}.

One has
\begin{multline*}(\cyl(\de\square_S[a_1,\dots,a_p]))_0 =
(\de\square_S[a_1,\dots,a_p])_0 \iso (\square_S[a_1,\dots,a_p])_0 \\ =
(\cyl(\square_S[a_1,\dots,a_p]))_0\end{multline*} and
\begin{multline*}(\cyl(\de\square_S[a_1,\dots,a_p]))_1 =
  (\de\square_S[a_1,\dots,a_p])_1 \p \{0,1\} \\ \iso
  (\square_S[a_1,\dots,a_p])_1 \p \{0,1\} =
  (\cyl(\square_S[a_1,\dots,a_p]))_1.\end{multline*} Therefore
$\overline{f}_0 = f_0$ and $\overline{f}_1 = f_1$ exist and are
unique.

The fact that $\ls{K}$ satisfies the HDA paradigm is going to be used
only now. Let $x\in (\cyl(\square_S[a_1,\dots,a_p]))_p$ with $p\geq
2$. Since
\[\cyl(\square_S[a_1,\dots,a_p]))_p = (\square_S[a_1,\dots,a_p]))_p \p \{0,1\}^p\]
one sees that all $\overline{f}(\de_i^\alpha x) = f(\de_i^\alpha x)$
are determined and the proof is complete.
\epf 

Hence the theorem: 

\bth \label{constr_model_HDA} There exists a unique combinatorial
model category structure on $\hda^\Sigma$ such that the class of
cofibrations is generated by $\I$ and such that a fibrant object is a
$\Lambda_{\square_S^{op}\set \ddownarrow
  !^S\Sigma}(\ls{!^S(\Sigma\p\{0,1\})},\varnothing,\I)$-injective
object of $\hda^\Sigma$.  All objects are cofibrant. \eth

Now we can state the last theorem: 

\bth \label{same2} There exists a Bousfield localization of
$\hda^\Sigma$ which is Quillen equivalent to the Bousfield
localization of the left determined model category structure of $\cts$
by the cubification functor.  \eth

\bpf By \cite[Theorem~9.5]{hdts}, the functor $\T:\square_S^{op}\set
\ddownarrow !^S\Sigma \to \whdts$ factors uniquely (up to isomorphism
of functors) as a composite
\[\square_S^{op}\set \ddownarrow
!^S\Sigma \stackrel{\sh_\Sigma} \longrightarrow \hda^{\Sigma}
\stackrel{\overline{\T}}\longrightarrow \whdts.\] Moreover, the
functor $\overline{\T}$ is a left adjoint. And it is not a left
Quillen adjoint with exactly the same proof as for
Proposition~\ref{not-left-Quillen}.  We work with the left Quillen
functor $\overline{\T}: \hda^\Sigma \to \bl_{\mathcal{S}}(\cts^+)$ like
in the proof of Theorem~\ref{same1}. We then just have to check that
the functor $\overline{\T}$ is, like the functor $\T$, homotopically
surjective.

Let $X$ be a cubical transition system. Let $n\geq 2$ and
$x_1,\dots,x_n\in\Sigma$. Then one has $(\square_S^{op}\set\ddownarrow
!^S\Sigma)(\square_n[x_1,\dots,x_n],\R(X)) \iso
\cts(\T(\square_n[x_1,\dots,x_n]),X)$ by the adjunction $\T\dashv \R$,
and
{\small
\begin{multline*}
  \cts(\T(\square_S[x_1,\dots,x_n]),X) =
  \cts(\T(\square_S[x_1,\dots,x_n] \sqcup_{\de
    \square_S[x_1,\dots,x_n]} \square_S[x_1,\dots,x_n]),X) \\ \iso
  \cts(\T(\square_S[x_1,\dots,x_n]),X)
  \sqcup_{\cts(\T(\de\square_S[x_1,\dots,x_n]),X)}\cts(\T(\square_S[x_1,\dots,x_n]),X)
\end{multline*}
}
by \cite[Proposition~9.3]{hdts} and because the contravariant functor
$\cts(-,X)$ is limit-preserv\-ing. So by the adjunction $\T\dashv \R$
again and because the contravariant functor
\[(\square_S^{op}\set\ddownarrow !^S\Sigma)(-,X)\] is limit-preserving,
one obtains the bijection
\begin{multline*}
(\square_S^{op}\set\ddownarrow
!^S\Sigma)(\square_n[x_1,\dots,x_n],\R(X)) \\\iso (\square_S^{op}\set\ddownarrow !^S\Sigma)(\square_S[x_1,\dots,x_n] \sqcup_{\de
  \square_S[x_1,\dots,x_n]} \square_S[x_1,\dots,x_n],\R(X)),
\end{multline*}
which means that $\R(X) \in \hda^\Sigma$. Therefore
$\overline{\T}(\R(X)) \iso \T(\R(X)) \iso \cub(X)$.
\epf

The reflection $\sh_\Sigma:\square_S^{op}\set\ddownarrow !^S\Sigma \to
\hda^\Sigma$ is a homotopically surjective left Quillen functor by
Theorem~\ref{restriction-model-structure}. So, by
\cite[Proposition~3.2]{MR1870516} there exists a set of maps
$\mathcal{Y}$ and a Quillen equivalence
\[\sh_\Sigma:\bl_{\mathcal{Y}}(\square_S^{op}\set\ddownarrow !^S\Sigma) \simeq
\hda^\Sigma.\] From Theorem~\ref{same2}, there exists a set of maps
$\mathcal{X}$ and a Quillen equivalence \[\overline{\T}:\bl_{\mathcal{X}}\hda^\Sigma
\simeq \bl_{\cub}(\cts).\] One obtains the Quillen equivalence
$\bl_{\mathcal{X}}\bl_{\mathcal{Y}}(\square_S^{op}\set\ddownarrow
!^S\Sigma) \simeq \bl_{\cub}(\cts)$.  Hence \[\T:\bl_{\mathcal{X}\cup
  \mathcal{Y}}(\square_S^{op}\set\ddownarrow !^S\Sigma) \simeq
\bl_{\cub}(\cts).\]

\end{document}